\documentclass[11pt]{article}

\usepackage[utf8]{inputenc}
\usepackage[T1]{fontenc}

\usepackage{amsmath,amsfonts,amssymb, amsthm}
\usepackage{hyperref}
\usepackage{xcolor}
\usepackage{float}
\usepackage{bm}
\usepackage{graphicx}
\usepackage{breakurl}
\usepackage{xifthen}
\usepackage{stmaryrd}
\usepackage{amsopn}
\usepackage{subcaption}

\usepackage{chngcntr}

\if@twoside
\else 
\fi

\theoremstyle{plain} 
\newtheorem{theorem}{Theorem}[section]

\theoremstyle{definition}

\theoremstyle{plain}
\newtheorem{remark}[theorem]{Remark}

\newcommand{\displaceP}{\overline{\bm{u}}}
\newcommand{\displace}{\bm{u}}
\newcommand{\displacet}{\bm{v}}

\newcommand{\displacePComp}{\tilde u}
\newcommand{\displaceComp}{u}
\newcommand{\displaceCompt}{v}
\newcommand{\displaceN}[1][]{
\ifthenelse{\isempty{#1}}
{u_3}
{u_{3 #1}}
}

\newcommand{\displaceNt}[1][]{
\ifthenelse{\isempty{#1}}
{v_3}
{v_{3 #1}}
}

\newcommand{\bendingStrainComp}{\kappa}
\newcommand{\membraneStrainComp}{\epsilon}
\newcommand{\bendingStrain}{\bm{\bendingStrainComp}}
\newcommand{\membraneStrain}{\bm{\membraneStrainComp}}

\newcommand{\bendingStrainFOrd}{\bm{\kappa}^{1}}

\newcommand{\bendingMoment}{\bm{M}}
\newcommand{\membraneForce}{\bm{N}}
\newcommand{\bendingMomentPComp}{\overline{M}}

\newcommand{\membraneForceComp}{N}
\newcommand{\bendingMomentt}{\bm{L}}
\newcommand{\membraneForcet}{\bm{K}}
\newcommand{\membraneForcetComp}{K}

\newcommand{\materialTensor}{\bm{C}}
\newcommand{\materialTensorComp}{C}
\newcommand{\materialTensorModBendingMoment}{\hat{\bm{C}}_{\bendingMoment}}
\newcommand{\materialTensorModMembraneForce}{\hat{\bm{C}}_{\membraneForce}}

\newcommand{\vp}{p}
\newcommand{\vpt}{q}
\newcommand{\vphi}{\bm{\phi}}
\newcommand{\vphit}{\bm{\psi}}
\newcommand{\vphitComp}{\psi}
\newcommand{\x}{\bm{x}}
\newcommand{\y}{\bm{y}}

\newcommand{\lagrange}{\bm{\lambda}}
\newcommand{\lagranget}{\bm{\mu}}

\newcommand{\lagrangetNDerh}{\mu_{n,h}}
\newcommand{\lagrangetTDerh}{\mu_{\tau,h}}

\newcommand{\map}{\bm{X}}
\newcommand{\mapMidSurf}{\bm{R}}

\newcommand{\mapMidSurfDef}{\bm{r}}

\newcommand{\baseVec}{\bm{A}}

\newcommand{\firstFfComp}{A}
\newcommand{\secFfComp}{B}

\newcommand{\JacDetMap}{\sqrt{a}} 

\newcommand{\Vprimal}{{\bm{W}}}

\newcommand{\VprimalN}{W_3}
\newcommand{\V}{{\bm{V}}}
\newcommand{\Q}{{\bm{Q}}}

\newcommand{\QN}{Q_3}

\newcommand{\VRegDecomp}{\bm{V}}
\newcommand{\SpLagrange}{\bm{\Lambda}}

\newcommand{\LtwobSym}{\bm{L}^2(\Omega)_{\mathrm{sym}}}

\newcommand{\grad}{\nabla}
\newcommand{\curl}{\operatorname{curl}}
\renewcommand{\div}{\operatorname{div}}

\newcommand{\hess}{\grad^2}

\newcommand{\bdiv}{\operatorname{Div}}
\newcommand{\bCurl}{\operatorname{Curl}}

\newcommand{\symbCurl}{\operatorname{symCurl}}

\newcommand{\laplace}{\Delta}

\newcommand{\BsplineSpace}{\mathcal{S}}

\DeclareMathOperator{\dotprod}{\cdot}
\newcommand{\pp}[2]{\partial_{#2} #1}
\newcommand{\dxi}{d\xi}

\newcommand{\Id}{\bm{I}}

\newcommand{\sectref}[1]{Section~\ref{#1}}
\newcommand{\appxref}[1]{Appendix}
\newcommand{\figref}[1]{Figure~\ref{#1}}
\newcommand{\tabref}[1]{Table~\ref{#1}}

\title{A new mixed isogeometric approach to Kirchhoff-Love shells}
\author{
	Katharina Rafetseder\footnote{Institute of Computational Mathematics, Johannes Kepler University Linz, 4040 Linz, Austria (rafetseder@numa.uni-linz.ac.at, zulehner@numa.uni-linz.ac.at).} \, and
	Walter Zulehner\footnotemark[1]
}

\begin{document}

\maketitle

\begin{abstract}
For Kichhoff-Love shell problems a new mixed formulation solely based on standard $H^1$ spaces is presented. This allows for flexibility in the construction of discretization spaces, e.g., standard $C^0$-coupling of multi-patch isogeometric spaces is sufficient. In terms of solution strategies, for iterative solvers efficient methods for standard second-order problems like multigrid can be used as building blocks of a preconditioner. Furthermore, a combination of the proposed mixed formulation of the bending part with a popular mixed formulation of the membrane part in order to avoid membrane locking is considered. The performance of both mixed formulations is demonstrated by numerical benchmark studies.

\smallskip
\noindent \textbf{Key words.} Kichhoff-Love shells, mixed methods, isogeometric shell elements, membrane locking

\end{abstract}


\section{Introduction}

In our previous work \cite{rafetseder_zulehner_2017} a new mixed variational formulation for the Kirchhoff plate bending problem with the bending moment tensor as additional unknown is derived. The new mixed formulation satisfies Brezzi's conditions and is equivalent to the original problem without additional convexity assumption on the domain. Furthermore, we obtain for polygonal domains an equivalent formulation of the Kirchhoff plate bending problem in terms of three (consecutively to solve) second-order elliptic problems. In \cite{rafetseder_zulehner_2017a} we extended the approach to domains whose boundaries are curvilinear polygons.  

The aim of this paper is to adapt the ideas of \cite{rafetseder_zulehner_2017} to the more complex situation of Kirchhoff-Love shells. 

Especially for shells, performing analysis directly on the geometry representation provided by the Computer Aided Design (CAD) model, and, thereby, avoiding unnecessary and costly geometry approximations is of essential importance. This is enabled by isogeometric analysis proposed by Hughes and co-workers in \cite{hughes_2005, hughes_2009}, where B-Splines or Non Uniform Rational B-Splines (NURBS) are used for the geometry description as well as for the representation of the unknown fields. 

Another feature of isogeometric discretizations is that they easily allow for inter-element continuity beyond the classical $C^0$-continuity within so-called patches. 
This permits the straight-forward construction of conforming isogeometric Kirchhoff-Love type shell elements, which are harder to obtain by means of standard Lagrange basis functions, since $C^1$-continuity is necessary as second derivatives appear. The contribution \cite{kiendl_2009} was one of the first to exploit this. 

For these approaches application on a single patch is straight-forward. However, for practical computations involving complex geometries usually geometry representations by means of several patches (multi-patches) are necessary. Then the continuity between patches is an issue, for which several techniques have been developed. 
In \cite{kiendl_2010} the so-called bending strip method is introduced, in which patches (called bending stripes) of fictitious material with bending stiffness only in direction transverse to the interface and zero membrane stiffness are added at patch interfaces. The crucial point is the choice of a reliable penalty parameter, the bending strip stiffness. Starting from this technique, alternative formulations removing the penalty parameter dependence have been proposed in \cite{goyal_2017}.
Alternatively, dG (discontinuous Galerkin) techniques can be used patch-wise, see, e.g., \cite{guo_ruess_2015}, where a variationally consistent Nitsche formulation, which weakly enforces coupling and continuity constraints among patches, is derived. Another approach are analysis suitable $C^1$ multi-patch isogeometric spaces, see, e.g., \cite{kapl_sangalli_takacs_2018}.

Note, all this techniques require special treatment of the patch interfaces. The essential contribution of this paper consists in a new mixed formulation of Kirchhoff-Love shells, with the bending moment tensor as new unknown, solely based on standard $H^1$ spaces. Therefore, classical $C^0$-coupling across patch interfaces is sufficient. In comparison with \cite{rafetseder_zulehner_2017}, the obtained mixed formulation differs in several aspects. In case of plates membrane and bending parts decouple, which is not the case for shells. Therefore, additionally the membrane part, which only involves first derivatives has to be taken into account. Furthermore, the bending strain includes in contrast to plates additional terms beside the Hessian of the transverse displacement involving geometry quantities and first derivatives of the displacement. In this paper it is shown how to overcome these difficulties and extend the techniques in \cite{rafetseder_zulehner_2017} to obtain again a formulation only based on $H^1$ spaces.

A build-in feature of the new mixed mixed formulation is the explicit computation of the bending moment tensor as separate unknown, which, however, leads to the drawback of having more unknowns. Beside the gained flexibility in the choice of discretization method, our mixed formulation also leads to new solution strategies. For the use of iterative solvers efficient methods for standard second-order problems like multigrid methods can be used as building blocks of a preconditioner.

It is a well known fact that also isogeometric Kirchhoff-Love shell elements exhibit significant membrane locking effects, see, e.g., \cite{echter_oesterle_bischoff_2013, echter_2013}. Therefore, we consider in this paper a combination of our new mixed formulation of the bending part with a popular method to avoid membrane locking by a mixed formulation of the membrane part, see \cite{echter_oesterle_bischoff_2013}. In the numerical experiments the combined method works well.

The paper is organized as follows: In \sectref{sec:primal_formulation} the Kirchhoff-Love shell formulation is introduced. \sectref{sec:mixed_formulation} contains a new mixed formulation and an extension of this approach to circumvent membrane locking. The new mixed formulation leads in a natural way to the construction of a new discretization method, which is introduced in \sectref{sec:discretization}. The paper closes with numerical experiments in \sectref{sec:numerical_experiments}.


\section{Kirchhoff-Love shell formulation}
\label{sec:primal_formulation}

The description of the Kirchhoff-Love shell model is based on the presentations in \cite{ciarlet_2005,  chapelle_bathe_2011, bischoff_wall_bletzinger_ramm_2004}.

\subsection{Differential geometry and shell kinematics}
Throughout this paper we use both index notation and absolute notation for the expression of vectors and tensors. Note, scalars are printed italic and vectors and tensors bold face italic. Quantities needed in the undeformed and deformed configuration are distinguished with capital and small letters, respectively.
Latin indices take values $\{1,2,3\}$ and Greek indices $\{1,2\}$. Superscripts indicate contravariant components and subscripts mark covariant components of vectors and tensors. Furthermore, Einstein's summation convention is applied to indices appearing twice within a product.

The 3D shell in the undeformed configuration is described by a mapping $\map$ from a parameter domain into the three-dimensional physical space, $\map: \Omega \times [-\frac{t}{2}, \frac{t}{2}] \rightarrow \mathbb{R}^3$, where $\Omega \subset \mathbb{R}^2$. The map $\map$ is defined in terms of curvilinear coordinates $\xi^i$ and is of the form
\begin{equation}
	\map(\xi^1, \xi^2, \xi^3) = \mapMidSurf(\xi^1, \xi^2) + \xi^3 \ \baseVec_3(\xi^1, \xi^2),
\end{equation}
where the mapping $\mapMidSurf: \Omega \rightarrow \mathbb{R}^3$ defines the midsurface of the shell, $\baseVec_3$ is the unit director, which is normal to the midsurface, and $t$ denotes the constant thickness . The map $\mapMidSurf$ will be assumed to be sufficiently smooth, e.g., $\mathcal{C}^3(\overline\Omega)$, cf. \cite{ciarlet_2005}.

The covariant base vectors are given by
\begin{equation*}
	\baseVec_\alpha = \pp{\mapMidSurf}{\alpha}, \qquad \baseVec_3 = \frac{\baseVec_1 \times \baseVec_2}{| \baseVec_1 \times \baseVec_2 |},
\end{equation*}
where here and in the following $\pp{}{\alpha} = \pp{}{\xi^\alpha}$. Contravariant base vectors are defined via the orthogonality condition
\begin{equation*}
	\baseVec_i \dotprod \baseVec^j = \delta^j_i.
\end{equation*}
Note that $\baseVec^3 = \baseVec_3$. The covariant and contravariant components $\firstFfComp_{\alpha \beta} $ and $\firstFfComp^{\alpha \beta}$ of the first fundamental form (metric tensor) of the shell surface, the Christoffel symbols $\Gamma^\sigma_{\alpha \beta}$, and the covariant and mixed components $\secFfComp_{\alpha \beta}$ and $\secFfComp^\beta_\alpha$ of the second fundamental form of the surface are then defined as follows:
\begin{equation}
\begin{alignedat}{4} 
	\firstFfComp_{\alpha \beta} &= \baseVec_\alpha \dotprod \baseVec_\beta, \quad 
	\firstFfComp^{\alpha \beta} = \baseVec^\alpha \dotprod \baseVec^\beta, \quad
	 \Gamma^\sigma_{\alpha \beta} = \baseVec^\sigma \dotprod \pp{\baseVec_\alpha}{\beta} \\
	\secFfComp_{\alpha \beta} &= \baseVec_3 \dotprod \pp{\baseVec_\alpha}{\beta}, \quad
	\secFfComp^\beta_\alpha = \firstFfComp^{\beta \sigma} \secFfComp_{\sigma \alpha}.
\end{alignedat}
\end{equation}
For later use let us define $\JacDetMap = \sqrt{\mathrm{det}(\firstFfComp_{\alpha \beta})}$.
We define analogously the above introduced quantities for the shell in the deformed configuration and denote them by the corresponding small letters.

The unknown we search for is the displacement $\displaceP = \mapMidSurfDef - \mapMidSurf$ at each point of the midsurface, which can be expressed in terms of covariant components $\displaceComp_i$ as follows
\begin{equation}
	\displaceP = \displaceComp_i \baseVec^i.
\end{equation}
We choose to solve for the vector of covariant components $\displace = (\displaceComp_i)$, since this simplifies the derivation of our new mixed formulation in \sectref{sec:Mmixed_formulation}. In the following we refer to $(\displaceComp_1, \displaceComp_2)$ and $\displaceComp_3$ as tangential and transverse part, respectively.

After linearization, we obtain expressions for the covariant components of the membrane strain tensor  
\begin{equation}
	\membraneStrainComp_{\alpha \beta}(\displace) = \frac{1}{2}(\displaceComp_{\alpha|\beta} + \displaceComp_{\beta|\alpha} ) - \secFfComp_{\alpha \beta} \ \displaceN
\end{equation}
and the bending strain tensor
\begin{equation}
	\bendingStrainComp_{\alpha \beta}(\displace) = \displaceN[|\alpha \beta] - \secFfComp^\sigma_\alpha \secFfComp_{\sigma \beta} \displaceN + \secFfComp^\sigma_\alpha \displaceComp_{\sigma|\beta} + \secFfComp^\tau_\beta \displaceComp_{\tau|\alpha} + \secFfComp^\tau_{\beta}|_\alpha \displaceComp_\tau.
\end{equation}
Here, $\displaceComp_{\alpha|\beta}$ denotes the covariant derivative 
$\displaceComp_{\alpha|\beta} = \pp{\displaceComp_\alpha}{\beta} - \Gamma^\sigma_{\alpha \beta} \displaceComp_\sigma$,
$\displaceN[|\alpha \beta]$ the second-order covariant derivative $\displaceN[|\alpha \beta] = \pp{\displaceN}{\alpha \beta} - \Gamma^\sigma_{\alpha\beta} \pp{\displaceN}{\sigma}$,
and $\secFfComp^\tau_\beta|_\alpha$ the covariant derivative of the second fundamental form
\begin{equation*}
	\secFfComp^\tau_\beta|_\alpha = \pp{\secFfComp^\tau_\beta}{\alpha} + \Gamma^\tau_{\alpha \sigma} \secFfComp^\sigma_\beta - \Gamma^\sigma_{\alpha \beta} \secFfComp^\tau_\sigma.
\end{equation*}
For a thorough derivation of this expressions see, e.g., \cite{ciarlet_2005}.

\subsection{Variational formulation}
The problem is posed in the parameter domain $\Omega$ with boundary $\Gamma$. In what follows let $(n_\alpha)$ and $\tau=(-n_2, n_1)$ represent the unit outer normal vector and the unit counterclockwise tangent vector to $\Gamma$, respectively. Furthermore, $\pp{}{n}$ denotes the normal derivative and $\pp{}{\tau}$ the tangential derivative along $\Gamma$. 
The shell is considered to be clamped on a part $\mapMidSurf(\Gamma_c)$, simply supported on $\mapMidSurf(\Gamma_s)$, and free on $\mapMidSurf(\Gamma_f)$, with $\Gamma = \Gamma_c \cup \Gamma_s \cup \Gamma_f$.
We define the displacement function space $\Vprimal$ by
\begin{equation}
\begin{aligned}
	\Vprimal = \{ \displacet = (\displaceCompt_i) \in H^1(\Omega) \times H^1(\Omega) \times H^2(\Omega): \ &\displaceCompt_i = 0, \ \pp{\displaceNt}{n} = 0 \ \text{on} \ \Gamma_c,\\
													   & \displaceCompt_1 = \displaceCompt_3 = 0 \ \text{on} \ \Gamma_s \}.
\end{aligned}
\end{equation}


The (primal) variational formulation is given as follows (cf. \cite{ciarlet_2005, chapelle_bathe_2011}) : find $\displace \in \Vprimal$ such that
\begin{equation}\label{eq:primal_formulation}
	\int_{\Omega} \left( t \ \membraneStrain(\displace) : \materialTensor : \membraneStrain(\displacet) + \frac{t^3}{12} \  \bendingStrain(\displace) : \materialTensor : \bendingStrain(\displacet) \right) \ \JacDetMap \ \dxi = \langle F, \displacet \rangle \quad \text{for all} \ \displacet \in \Vprimal,
\end{equation}
where $\bm{A}:\bm{B}$ denotes the double contraction of two tensors $\bm{A}$ and $\bm{B}$ and $\dxi = d\xi^1 d \xi^2$.
Here, the right-hand side is given by $\langle F, \displacet \rangle = \int_\Omega \bm{f} \dotprod \displacet \ \JacDetMap \ \dxi$ and the contravariant components of the fourth-order material tensor read
\begin{equation*}
	\materialTensorComp^{\alpha \beta \sigma \tau} = \frac{E}{2(1+\nu)} \left(\firstFfComp^{\alpha \sigma} \firstFfComp^{\beta \tau} + \firstFfComp^{\alpha \tau} \firstFfComp^{\beta \sigma} + \frac{2 \nu}{1-\nu} \firstFfComp^{\alpha \beta} \firstFfComp^{\sigma \tau}\right),
\end{equation*}
where $E$ and $\nu$ are Young's modulus and Poisson's ratio of the material, respectively.

Following, e.g., \cite{adams_fournier_2003}, here and throughout the paper  $L^2(\Omega)$ and $H^m(\Omega)$ denote the standard Lebesgue and Sobolev spaces of functions on $\Omega$ with corresponding norms $\|.\|_{0}$ and $\|.\|_{m}$ for positive integers $m$. Moreover, $\LtwobSym$ denotes the space of symmetric second-order tensors given by
\begin{equation*}
	\LtwobSym = \{ \membraneForcet : \membraneForcetComp^{\alpha\beta} = \membraneForcetComp^{\beta\alpha} \in L^2(\Omega) \}.
\end{equation*}
For scalars $\displaceCompt$, vectors $\vphit$, and second-order tensors $\membraneForcet$ the first order differential operators with respect to the tangential coordinates $\xi^\alpha$ are defined as follows:
\begin{align*}
	\grad \displaceCompt &= 
	\begin{pmatrix}
	\pp{\displaceCompt}{1} \\
	\pp{\displaceCompt}{2}
	\end{pmatrix},\quad
	&\curl \displaceCompt &=
	\begin{pmatrix}
	\pp{\displaceCompt}{2} \\
	-\pp{\displaceCompt}{1}
	\end{pmatrix},\\
	\grad \vphit &=
	\begin{pmatrix}
	\pp{\vphitComp_{1}}{1}  & \pp{\vphitComp_{1}}{2}  \\
	\pp{\vphitComp_{2}}{1} & \pp{\vphitComp_{2}}{2}  \\
	\end{pmatrix},\quad
	&\bCurl \vphit&=
	\begin{pmatrix}
	\pp{\vphitComp_{1}}{2} & -\pp{\vphitComp_{1}}{1}  \\
	\pp{\vphitComp_{2}}{2} & -\pp{\vphitComp_{2}}{1}  \\
	\end{pmatrix},\\
	\div \vphit &= \pp{\vphitComp_1}{1} + \pp{\vphitComp_2}{2}, \quad 
	&\bdiv \membraneForcet &= 
	\begin{pmatrix}
	\pp{\membraneForcetComp_{11}}{1} + \pp{\membraneForcetComp_{12}}{2}  \\
	\pp{\membraneForcetComp_{21}}{1} + \pp{\membraneForcetComp_{22}}{2}  \\
	\end{pmatrix}.
\end{align*}
Moreover, the symmetric $\bCurl$ is introduced by
\begin{equation*}
	\symbCurl \vphit = \frac{1}{2}( \bCurl \vphit + (\bCurl \vphit)^T).
\end{equation*}


\section{Two mixed variational formulations}
\label{sec:mixed_formulation}
In this section we derive two mixed formulations by introducing stress resultants as new unknowns.

\subsection{\texorpdfstring{$\bendingMoment$}{M}-mixed formulation}
\label{sec:Mmixed_formulation}
We introduce as new unknown the bending moment tensor $\bendingMoment$, which is related to the bending strain through the constitutive equation
\begin{equation*}
	\bendingMoment = \JacDetMap \frac{t^3}{12} \ \materialTensor : \bendingStrain(\displace) = \materialTensorModBendingMoment : \bendingStrain(\displace), \ \text{with} \ \materialTensorModBendingMoment = \JacDetMap \frac{t^3}{12} \ \materialTensor.
\end{equation*}
Note, in contrast to standard notation we additionally include in $\bendingMoment$ the geometry measure $\JacDetMap$. This leads to the preliminary $\bendingMoment$-mixed formulation: find $\bendingMoment \in \LtwobSym$ and $\displace \in \Vprimal$ such that
\begin{equation}\label{eq:first_Mmixed_formulation}
\begin{alignedat}{4} 
	& \int_\Omega (\materialTensorModBendingMoment^{-1} : \bendingMoment) : \bendingMomentt \ \dxi & & - \int_\Omega \bendingStrain(\displace) : \bendingMomentt \ \dxi & & = 0\\
	& -\int_{\Omega} \bendingMoment : \bendingStrain(\displacet)  \ \dxi & & - c(\displace, \displacet) & &  = -\langle F, \displacet \rangle
\end{alignedat}
\end{equation}
for all $\bendingMomentt \in \LtwobSym$ and $\displacet \in \Vprimal$, with $c(\displace, \displacet) $ given by the membrane part
\begin{equation*}
	c(\displace, \displacet) = \int_{\Omega} t \  \membraneStrain(\displace) : \materialTensor : \membraneStrain(\displacet) \ \JacDetMap \ \dxi.
\end{equation*}
In order to make equations more compact, we use in the following at some points the short notation $(.,.)$ for the $L^2$-inner product on $\Omega$ instead of an integral and put the material tensor as subscript, i.e., $(\materialTensorModBendingMoment^{-1} \bendingMoment, \bendingMomentt) = (\bendingMoment, \bendingMomentt)_{\materialTensorModBendingMoment^{-1}}$.

The goal of the remainder of this section is to derive a reformulation of \eqref{eq:first_Mmixed_formulation} allowing us to replace the displacement space $\Vprimal$ by a space that uses $H^1(\Omega)$ for all three components of the displacement. In our previous work \cite{rafetseder_zulehner_2017} a new mixed formulation for Kirchhoff plates using $H^1(\Omega)$ for the vertical deflection is introduced. An extension of this approach to shells is possible, since the only term involving second-order derivatives is the Hessian of the transverse displacement $\displaceN$.

We divide the bending strain $\bendingStrain(\displacet)$ into the Hessian of the transverse displacement $\displaceNt$ and the remaining terms that only involve first-order derivatives of $\displacet$ denoted by $\bendingStrainFOrd(\displacet)$, i.e.,
\begin{equation*}
	\bendingStrain(\displacet) = \hess \displaceNt + \bendingStrainFOrd(\displacet).
\end{equation*}
With this notation we can rewrite the second line in \eqref{eq:first_Mmixed_formulation} separating the integral involving $\hess \displaceNt$ and put all remaining terms into the right-hand side
$\langle G(\bendingMoment,c,F), \displacet \rangle = (\bendingMoment, \bendingStrainFOrd(\displacet)) + c(\displace, \displacet) -\langle F, \displacet \rangle$ leading to
\begin{equation*}
	-\int_\Omega \bendingMoment : \hess \displaceNt \ \dxi = \langle G(\bendingMoment,c,F) , \displacet \rangle,
\end{equation*}
or in strong form 
\begin{equation}\label{eq:secLine_first_Mmixed_formulation}
	-\div\bdiv \bendingMoment = G(\bendingMoment,c,F).
\end{equation}
The main idea is the following ansatz for $\bendingMoment$:
\begin{equation*}
	\bendingMoment = \vp\Id + \bendingMoment_0,
\end{equation*}
where $\div\bdiv \bendingMoment_0 = 0$ and $\Id$ is the identity matrix. Plugging in \eqref{eq:secLine_first_Mmixed_formulation} the just stated representation provides
\begin{equation}\label{eq:thirdLine_new_Mmixed_formulation}
	-\div\bdiv (\vp\Id) = G(\bendingMoment,c,F), \ \text{or equivalently} \ - \laplace \vp = G(\bendingMoment,c,F).
\end{equation}
Therefore, it is sufficient to consider $\vp \in H^1(\Omega)$ for the corresponding weak form. 
The second essential ingredient is a characterization of the elements in the kernel of $\div\bdiv$. According to \cite{beirao_niiranen_stenberg_2007, huang_huang_xu_2011, krendl_rafetseder_zulehner_2016} there is a potential function $\vphi \in (H^1(\Omega))^2$ such that
\begin{equation*}
	\bendingMoment_0 = \symbCurl \vphi.
\end{equation*}
\begin{remark}
	The first step can be viewed as homogenization. In literature, another form of homogenization has been considered. The classical Helmholtz decomposition $\bendingMoment = \grad^2 \vp + \symbCurl\vphi$, see, e.g., \cite{beirao_niiranen_stenberg_2007, huang_huang_xu_2011}, has the same second component. However, the first component is different and requires the solution of a fourth-order problem, which brings no benefit. In contrast, the decomposition introduced here only requires the solution of a second-order Poisson problem for the first component.
	
	Note the analogy to the well-known characterization of the stress tensor $\bm{\sigma}$ with $\bdiv \bm{\sigma} = 0$ by means of the Airy stress function $\varphi$ in $2$D
	\begin{equation*}
		\bm{\sigma} = \bCurl\curl \varphi,
	\end{equation*}
	and a similar result in $3$D with the Beltrami stress functions. 
\end{remark}

Summing up, we have the following representation of $\bendingMoment$:
\begin{equation*}
	\bendingMoment = \vp \Id + \symbCurl \vphi.
\end{equation*}
With this representation the preliminary $(\bendingMoment, \displace)$ problem in \eqref{eq:first_Mmixed_formulation} becomes a formulation in $(\vp, \vphi, \displace)$, later referred to as $\bendingMoment$-mixed formulation, given as follows: find $(\vp, \vphi)\in \VRegDecomp$ and $\displace \in \Q$ such that
\begin{small}
\begin{equation}\label{eq:new_Mmixed_formulation} 
\begin{alignedat}{5}
	&(\vp \Id, \vpt \Id)_{\materialTensorModBendingMoment^{-1}} && +(\symbCurl \vphi, \vpt \Id)_{\materialTensorModBendingMoment^{-1}} && +(\grad \displaceN, \grad \vpt) - (\bendingStrainFOrd(\displace), \vpt \Id)   &&= 0\\ 
	&(\vp \Id,\symbCurl \vphit)_{\materialTensorModBendingMoment^{-1}} && +(\symbCurl \vphi, \symbCurl \vphit)_{\materialTensorModBendingMoment^{-1}}  && - (\bendingStrainFOrd(\displace), \symbCurl \vphit)   &&= 0\\
	&(\grad \vp, \grad \displaceNt) - ( \vp \Id, \bendingStrainFOrd(\displacet)) &&  - (\symbCurl \vphi, \bendingStrainFOrd(\displacet))  && -c(\displace, \displacet) && = -\langle F, \displacet \rangle.
\end{alignedat}
\end{equation}
\end{small}
for all $(\vpt, \vphit)\in \VRegDecomp$ and $\displacet \in \Q$. 
Here, the third line is given by the weak form of \eqref{eq:thirdLine_new_Mmixed_formulation}, where we plug in the right-hand side $G(\bendingMoment,c,F)$ the representation of $\bendingMoment$. The first two lines follow from the first line in \eqref{eq:first_Mmixed_formulation} using an analogous representation of the test functions $\bendingMomentt = \vpt \Id + \symbCurl \vphit$ and splitting of the bending strain $\bendingStrain(\displace) = \hess \displaceN + \bendingStrainFOrd(\displace)$.

The original displacement space $\Vprimal$ is replaces by the space $\Q$ defined by
\begin{equation}
\begin{aligned}
	\Q &= \{ \displacet = (\displaceCompt_i) \in H^1(\Omega) \times H^1(\Omega) \times H^1(\Omega): \ \displaceCompt_i = 0 \ \text{on} \ \Gamma_c, \displaceCompt_1 = \displaceCompt_3 = 0 \ \text{on} \ \Gamma_s \}.
\end{aligned}
\end{equation}
The definition of the right boundary conditions for $\vp$ and $\vphi$ is a subtle issue. It turns out that the space $\VRegDecomp$ is given by the subset of $ (\vpt, \vphit) \in \QN \times (H^1(\Omega))^2$, with
\begin{equation*}
	\QN = \{ \displaceNt \in H^1(\Omega): \displaceCompt_3 = 0 \ \text{on} \ \Gamma_c \cup \Gamma_s \},
\end{equation*}
satisfying the boundary condition
\begin{equation}\label{eq:coupling_condition}
	\langle \pp{\vphit} {\tau} ,\grad \displaceNt \rangle_\Gamma + \int_\Gamma \vpt \ \pp{\displaceNt}{n} \ ds = 0 \quad \text{for all} \ \displaceNt\in \VprimalN,
\end{equation}
where $\langle ., . \rangle_\Gamma$ denotes the duality product on $\Gamma$ and
\begin{equation*}
	\VprimalN = \{ \displaceNt \in H^2(\Omega): \ \displaceCompt_3 = 0, \ \pp{\displaceCompt_3}{n} = 0 \ \text{on} \ \Gamma_c, \quad \displaceCompt_3 = 0 \ \text{on} \ \Gamma_s \}.
\end{equation*}
By \eqref{eq:coupling_condition} the functions $\vpt$ and $\vphit$ are coupled. Therefore, we refer to \eqref{eq:coupling_condition} in the following as coupling condition.
In case $\vphit$ is sufficiently smooth, e.g., $\pp{\vphit}{\tau} \in L^2(\Gamma)$, the coupling condition can be rewritten as
\begin{equation}\label{eq:coupling_condition_smooth} 
	\int_{\Gamma_s \cup \Gamma_f} (\pp{\vphit}{\tau}\dotprod n + \vpt) \pp{\displaceNt}{n} \ ds + \int_{\Gamma_f} \pp{\vphit}{\tau}\dotprod \tau \ \pp{\displaceNt}{\tau} \ ds = 0 \quad \text{for all} \ \displaceNt \in \VprimalN,
\end{equation}
where we use the representation $\grad \displaceNt = (\pp{\displaceNt}{n}) \, n + (\pp{\displaceNt}{\tau}) \, \tau$ and incorporate the boundary conditions for $\displaceNt \in \VprimalN$. This condition reads in explicit form
\begin{equation}\label{eq:boundaryCond_p_phi}
\begin{aligned}
	&\pp{\vphit} {\tau} \dotprod n = -\vpt \quad && \text{on} \ \Gamma_s\\
	&\pp{^2\vphit} {\tau} \dotprod \tau  = 0,  \quad \pp{\vphit} {\tau} \dotprod n = -\vpt \quad && \text{on} \ \Gamma_f.
\end{aligned}	
\end{equation}
\begin{remark}
	Note, in $\Q$ compared with $\Vprimal$ only those boundary conditions for the transverse displacement $\displaceN$ which are also available in $H^1(\Omega)$ are prescribed. In the formulation \eqref{eq:new_Mmixed_formulation} the originally essential boundary condition $\pp{\displaceN}{n} = 0$ becomes a natural condition. In the primal formulation \eqref{eq:primal_formulation} only natural boundary conditions are imposed for the bending moment tensor $\bendingMoment$. Those natural boundary conditions corresponding to the normal-normal component of $\bendingMoment$ and the corner conditions for the normal-tangential component of $\bendingMoment$ become the essential conditions \eqref{eq:boundaryCond_p_phi} for $(\vp,\vphi)$ and the remaining condition remains natural. For further information we refer the reader to \cite{rafetseder_zulehner_2017}.
\end{remark}
\begin{remark}
	The mixed formulation \eqref{eq:new_Mmixed_formulation} is different from the formulation obtained in \cite{rafetseder_zulehner_2017}. In case of plates membrane and bending parts decouple and for the bending part a decomposition into three consecutively to solve second-order problems is obtained. This is no longer possible for shells.
\end{remark}

Problem \eqref{eq:new_Mmixed_formulation} has the typical structure of a saddle point problem: find $\x = (\vp, \vphi) \in \VRegDecomp$ and $\displace \in \Q$ such that
\begin{equation*}
\begin{alignedat}{4}
  & a(\x, \y) & & + b(\y,\displace) & & = 0  & \quad & \text{for all} \ \y = (\vpt, \vphit) \in \VRegDecomp , \\
  & b(\x, \displacet) & & - c(\displace, \displacet) & & = - \langle F , \displacet \rangle & \quad & \text{for all} \ \displacet \in \Q.
\end{alignedat}
\end{equation*}
As in \cite{rafetseder_zulehner_2017}, for the new mixed formulation \eqref{eq:new_Mmixed_formulation} the following result providing well-posedness and equivalence to the primal formulation holds: 
\begin{theorem}
	Let $\Omega$ be simply connected. The mixed formulation \eqref{eq:new_Mmixed_formulation} is well-posed, i.e., existence and uniqueness of a solution $((\vp,\vphi),\displace) \in \VRegDecomp \times \Q$ is guaranteed. Moreover, $\displace$ is the solution of the primal problem \eqref{eq:primal_formulation}, with the bending moment tensor $\bendingMoment = \materialTensorModBendingMoment : \bendingStrain(\displace)$ given by $\bendingMoment = \vp\Id + \symbCurl \vphi$.
\end{theorem}
The proof follows along the lines of \cite[Theorem 3.6]{rafetseder_zulehner_2017} and is omitted here. For more details on the mathematical foundation of the new mixed formulation we refer to the thorough derivation of an analogous mixed formulation for Kirchhoff plates in \cite{rafetseder_zulehner_2017}.

\subsection{\texorpdfstring{$\bendingMoment$-$\membraneForce$}{MN}-mixed formulation}
\label{sec:MNmixed_formulation}
In order to alleviate membrane locking one popular concept (among many others) is to consider a mixed formulation with the membrane force tensor $\membraneForce$ as new unknown, see, e.g., \cite{echter_oesterle_bischoff_2013,chapelle_stenberg_1999}. An adoption of this approach allows for a well matching extension of the just introduced formulation.

We additionally introduce as new unknown the membrane force tensor $\membraneForce$, which is connected to the membrane strain through the constitutive equation
\begin{equation*}
	\membraneForce = \JacDetMap \ t \ \materialTensor : \membraneStrain(\displace) = \materialTensorModMembraneForce : \membraneStrain(\displace), \ \text{with} \ \materialTensorModMembraneForce = \JacDetMap \ t \ \materialTensor.
\end{equation*}
Note, again we include in $\membraneForce$ the geometry measure $\JacDetMap$. This leads to the preliminary $\bendingMoment$-$\membraneForce$-mixed formulation: find $\bendingMoment \in \LtwobSym$, $\membraneForce \in \LtwobSym$ and $\displace \in \Vprimal$ such that
\begin{equation}\label{eq:first_NMmixed_formulation}
\begin{alignedat}{5} 
	& \int_\Omega (\materialTensorModBendingMoment^{-1} : \bendingMoment) : \bendingMomentt \ \dxi & &  & & - \int_\Omega \bendingStrain(\displace) : \bendingMomentt \ \dxi & & = 0\\
	&  && \int_\Omega (\materialTensorModMembraneForce^{-1} : \membraneForce) : \membraneForcet \ \dxi & & - \int_\Omega \membraneStrain(\displace) : \membraneForcet \ \dxi & & = 0\\
	& -\int_{\Omega} \bendingMoment : \bendingStrain(\displacet)  \ \dxi & & -\int_{\Omega} \membraneForce : \membraneStrain(\displacet)  \ \dxi & & & &  = -\langle F, \displacet \rangle
\end{alignedat}
\end{equation}
for all $\bendingMomentt \in \LtwobSym$, $\membraneForcet \in \LtwobSym$ and $\displacet \in \Vprimal$. Analogously as above, we can reformulate the preliminary ($\bendingMoment$, $\membraneForce$, $\displace$) problem  \eqref{eq:first_NMmixed_formulation} in terms of ($\vp$, $\vphi$, $\membraneForce$, $\displace$), later referred to as $\bendingMoment$-$\membraneForce$-mixed formulation: find $(\vp, \vphi)\in \VRegDecomp$, $\membraneForce\in \LtwobSym$ and $\displace \in \Q$ such that
\begin{small}
\begin{alignat*}{5}
	&(\vp \Id, \vpt \Id)_{\materialTensorModBendingMoment^{-1}} && +(\symbCurl \vphi, \vpt \Id)_{\materialTensorModBendingMoment^{-1}} &&  && +(\grad \displaceN, \grad \vpt) - (\bendingStrainFOrd(\displace), \vpt \Id)   &&= 0\\
	&(\vp \Id,\symbCurl \vphit)_{\materialTensorModBendingMoment^{-1}} && +(\symbCurl \vphi, \symbCurl \vphit)_{\materialTensorModBendingMoment^{-1}} &&  && - (\bendingStrainFOrd(\displace), \symbCurl \vphit)   &&= 0\\
	& && && \ (\membraneForce, \membraneForcet)_{\materialTensorModMembraneForce^{-1}} && -(\membraneStrain(\displace), \membraneForcet) && = 0\\
	&(\grad \vp, \grad \displaceNt) - ( \vp \Id, \bendingStrainFOrd(\displacet)) &&  - (\symbCurl \vphi, \bendingStrainFOrd(\displacet))  && -(\membraneForce, \membraneStrain(\displacet)) && && = -\langle F, \displacet \rangle.
\end{alignat*}
\end{small}
for all $(\vpt, \vphit)\in \VRegDecomp$, $\membraneForcet \in \LtwobSym$ and $\displacet \in \Q$


\section{The discretization method}
\label{sec:discretization}

In this section we first construct a conforming discretization space for the displacement $\displace$, i.e.,
\begin{equation*}
	\Q_h \subset \Q \subset (H^1(\Omega))^3.\\
\end{equation*}
Note, only $C^0$ basis functions are required, so the continuity requirements are easily satisfied with standard basis functions. In the following we consider isogeometric B-spline discretization spaces with degree $p\geq1$. For $p=1$, the discretization space coincides with the standard isoparametric finite element space of continuous and piecewise bilinear elements. We define discretization spaces on patch level and for our formulation continuity between patches does not require extra attention, since standard $C^0$-coupling is sufficient. 

We denote by $\BsplineSpace^{p_1,p_2}_{\alpha_1,\alpha_2}$ the tensor product B-spline space defined as
\begin{equation}\label{eq:BsplineSpace}
	\BsplineSpace^{p_1,p_2}_{\alpha_1,\alpha_2} = \BsplineSpace^{p_1}_{\alpha_1} \otimes \BsplineSpace^{p_2}_{\alpha_2},
\end{equation}
where $\BsplineSpace^{p}_{\alpha}$ is the one-dimensional B-spline space with degree $p$ and $\alpha$ continuous derivatives across interior knots; see, e.g, \cite{cotrell_hughes_brazilevs_2009,daVeiga_buffa_sangalli_vazquez_2014} for further information.

We use equal order discretization spaces for the three components of the displacement $\displace$ and incorporate the essential boundary conditions, which brings us to the definition:
\begin{equation*}
	\Q_h = (\BsplineSpace^{p,p}_{\alpha,\alpha})^3 \cap \Q,
\end{equation*}
where $\alpha = p-1$, i.e., maximum smoothness at interior knots. 

For $\VRegDecomp$, the space of the auxiliary variables $(\vp, \vphi)$, the construction of a conforming discretization space is more involved, since the coupling condition \eqref{eq:coupling_condition} has to be taken into account. Therefore, we first disregard the coupling condition and construct a conforming discretization space of $\hat\VRegDecomp = \QN \times (H^1(\Omega))^2$. Using equal order discretization spaces for $\vp$ and $\vphi$ we receive 
\begin{equation*}
	\hat\VRegDecomp_h = (\BsplineSpace^{p,p}_{\alpha,\alpha} \times (\BsplineSpace^{p,p}_{\alpha,\alpha})^2) \cap \hat\VRegDecomp.
\end{equation*}

The space $\VRegDecomp_h$ is defined as the subset of $\y_h = (\vpt_h, \vphit_h) \in \hat\VRegDecomp_h$ satisfying a discrete version of the coupling condition \eqref{eq:coupling_condition_smooth}
\begin{equation*}
	\VRegDecomp_h = \{ \y_h = (\vpt_h, \vphit_h)\in\hat\VRegDecomp_h : d(\y_h, \lagranget_h) = 0 \quad \text{for all} \ \lagranget_h=(\lagrangetTDerh,\lagrangetNDerh) \in \SpLagrange_h \},
\end{equation*}
with
\begin{equation}\label{eq:coupling_condition_discrete}
	d(\y_h, \lagranget_h) = \int_{\Gamma_s \cup \Gamma_f} (\pp{\vphit_h}{\tau}\dotprod n + \vpt_h) \ \lagrangetNDerh \ ds + \int_{\Gamma_f} \pp{\vphit_h}{\tau}\dotprod \tau \ \lagrangetTDerh \ ds.
\end{equation}
The test functions $(\lagrangetTDerh,\lagrangetNDerh)$ are discrete representations of $(\pp{\displaceNt}{\tau}, \pp{\displaceNt}{n})$ for $\displaceNt\in\VprimalN$ at $\Gamma$. Therefore, $\SpLagrange_h$ is chosen as the space of restrictions of functions from $(\BsplineSpace^{p-1,p-1}_{\alpha,\alpha})^2$ to $\Gamma$, where at corner points of the boundary $\lagrangetTDerh$ and $\lagrangetNDerh$ have to be coupled appropriately, see \cite{rafetseder_zulehner_2017a} for details. Then the discrete version of \eqref{eq:new_Mmixed_formulation} reads:
find $\x_h = (\vp_h, \vphi_h) \in \VRegDecomp_h$ and $\displace_h \in \Q_h$ such that
\begin{equation*}
\begin{alignedat}{4}
  & a(\x_h, \y_h) & & + b(\y_h,\displace_h) & & = 0  & \quad & \text{for all} \ \y_h = (\vpt_h, \vphit_h) \in \VRegDecomp_h, \\
  & b(\x_h, \displacet_h) & & - c(\displace_h, \displacet_h) & & = - \langle F , \displacet_h \rangle & \quad & \text{for all} \ \displacet_h \in \Q_h.
\end{alignedat}
\end{equation*}

We do not explicitly build in the coupling condition by constructing a basis of the space $\VRegDecomp_h$, but incorporate it implicitly, by replacing $\V_h$ by $\hat\V_h$ and adding \eqref{eq:coupling_condition_discrete} as additional constraint. Since $\hat\V_h$, $\Q_h$ and $\SpLagrange_h$ are finite dimensional, the bilinear forms $a$, $b$, $c$ and $d$ can be represented as matrices $\bm{A}_h$, $\bm{B}_h$, $\bm{C}_h$ and $\bm{D}_h$ acting on vectors of real numbers $\underline{\x}_h$, $\underline{\displace}_h$ and $\underline{\lagrange}_h$ representing the elements in $\hat\V_h$, $\Q_h$ and $\SpLagrange_h$, respectively, with respect to the chosen basis. In this matrix-vector notation the resulting system reads
\begin{equation*}
\begin{alignedat}{4}
	&\bm{A}_h \underline{\x}_h &&+ \bm{B}_h^T \underline{\displace}_h &&+ \bm{D}_h^T \underline{\lagrange}_h &&= 0,\\
	&\bm{B}_h \underline{\x}_h &&- \bm{C}_h   \underline{\displace}_h && &&= \bm{f}_h,\\
	&\bm{D}_h \underline{\x}_h && && &&=0,
\end{alignedat}
\end{equation*}
where $\bm{B}_h^T$ and $\bm{D}_h^T$ denote the transposed matrices.


For the second, the $\bendingMoment$-$\membraneForce$-mixed formulation we consider for $\displace$ and $(\vp, \vphi)$ the discretization introduced above. 
For the additional unknown, the membrane force $\membraneForce$, we use the discretization space proposed in the Hybrid Stress (HS) method presented in \cite{echter_oesterle_bischoff_2013}. The basis for the contravariant components of $\membraneForce$ is given by
\begin{equation*}
	{\membraneForceComp}^{11}_h \in \BsplineSpace^{p-1,p}_{\alpha-1,\alpha}, \quad
	{\membraneForceComp}^{22}_h \in \BsplineSpace^{p,p-1}_{\alpha,\alpha-1}, \quad
	{\membraneForceComp}^{12}_h \in \BsplineSpace^{p-1,p-1}_{\alpha-1,\alpha-1}.
\end{equation*}


\section{Numerical experiments}
\label{sec:numerical_experiments}
In the first part of this section we demonstrate that the $\bendingMoment$-mixed formulation (with the $H^1(\Omega)$ conforming discretizations proposed in the previous section) works by testing it with the three benchmark problems of the well-known shell obstacle course \cite{belytschko_1985}, consisting of two cylindrical shells and one spherical shell.
In the second part we show that the $\bendingMoment$-$\membraneForce$-mixed formulation introduced in \sectref{sec:MNmixed_formulation} works well.

For all problems in this paper the undeformed midsurface is modeled exactly by non-uniform rational B-splines (NURBS). For the discretization of the unknowns standard B-spline spaces are used, see \sectref{sec:discretization}. Mesh density is characterized by the number of control points per edge. In case several congruent patches are used to describe the surface the corresponding number for one patch is used. 



The implementation is done in the framework of G+Smo ("Geometry + Simulation Modules"), an object-oriented C++ library, see \url{https://ricamsvn.ricam.oeaw.ac.at/trac/gismo/wiki/WikiStart}. In all experiments a sparse direct solver is used.

\subsection{\texorpdfstring{$\bendingMoment$}{M}-mixed formulation}

\subsubsection{Scordelis-Lo roof}
\label{sec:roof}
\begin{figure}[!ht]
	\centering
	\begin{subfigure}[b]{0.44\textwidth}
		\includegraphics[width=\textwidth]{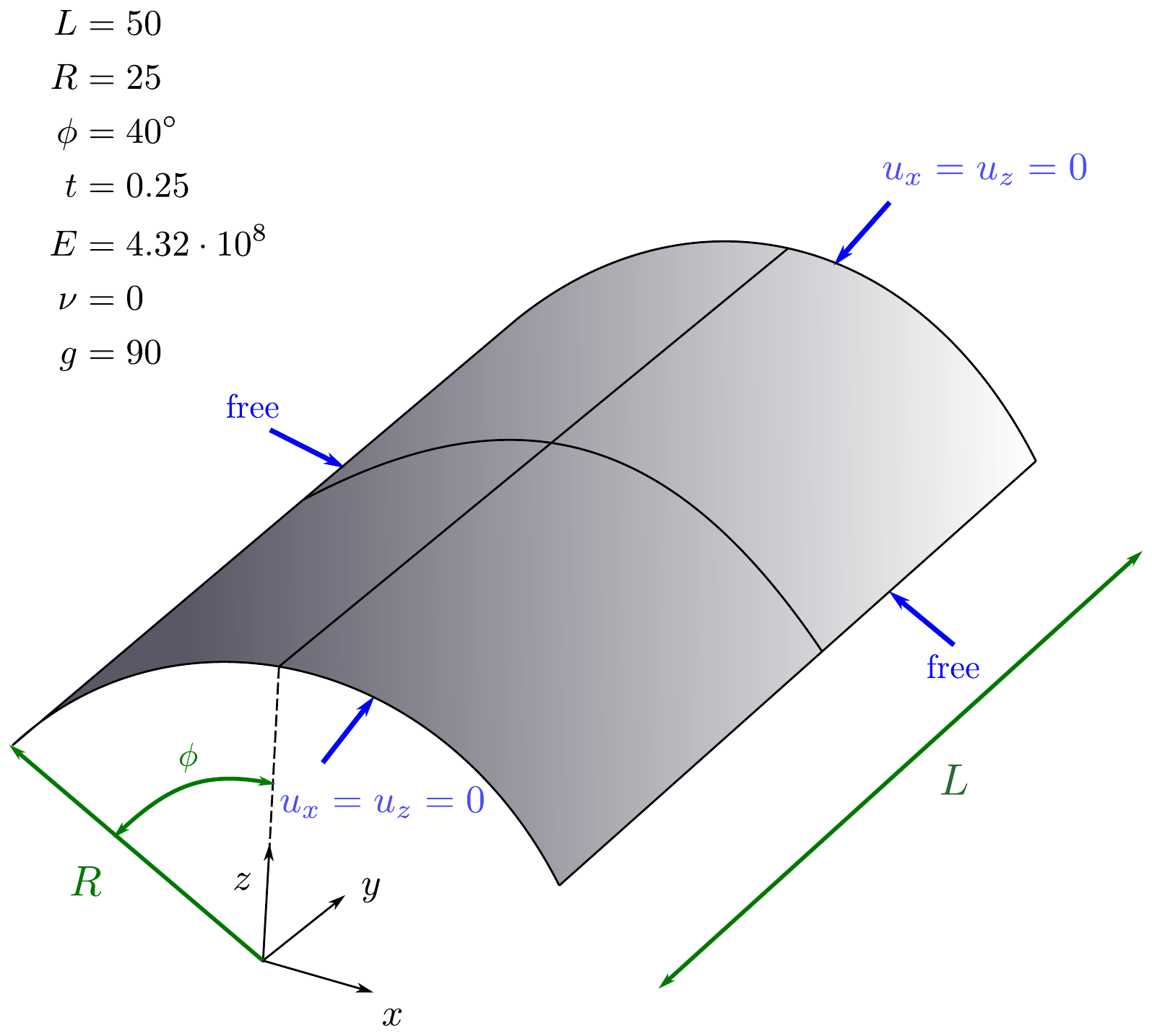}
		\caption{Problem setup.}
		\label{fig:roof_geometry}
	\end{subfigure}
	\begin{subfigure}[b]{0.55\textwidth}
		\includegraphics[width=\textwidth]{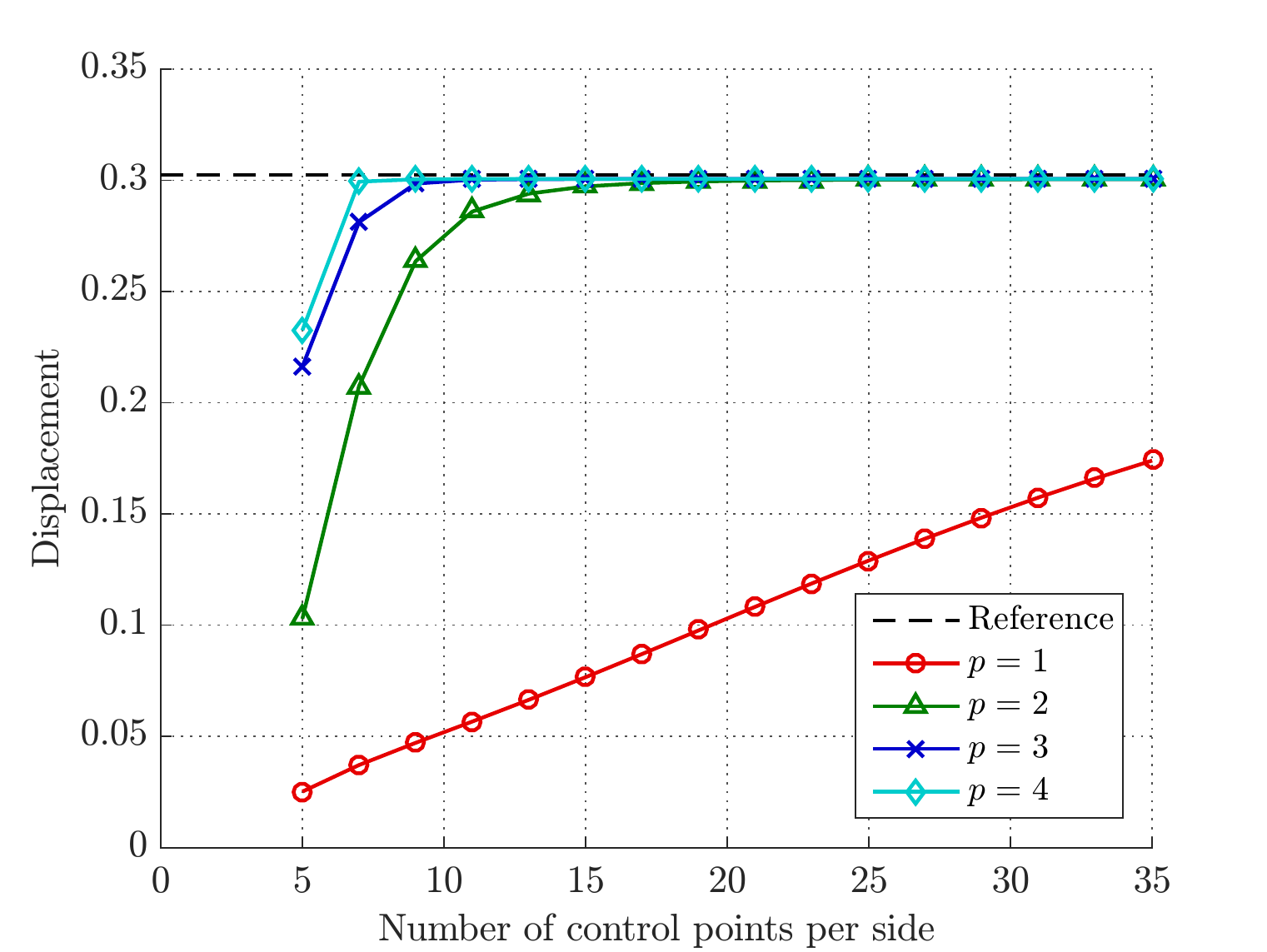}
		\caption{Displacement convergence of $\bendingMoment$-mixed formulation.}
		\label{fig:roof_conv_Mmixed}
	\end{subfigure}
	\caption{Scordelis-Lo roof.}
\end{figure}

The problem setup of the Scordelis-Lo roof benchmark is shown in \figref{fig:roof_geometry}. The structure is supported with rigid diaphragms at both ends and the side edges are free. This setup is realized by imposing homogeneous boundary conditions for the displacement of the form $\displacePComp_x = \displacePComp_z = 0$, leading to the conditions $\displaceComp_1 = \displaceComp_3 = 0$ for the covariant components of the displacement. The roof has moderate slenderness $\frac{R}{t} = 100$, with radius $R$ and thickness $t$ as defined in \figref{fig:roof_geometry}, and is subject to a uniform vertical load of $g = 90$ per unit area. This configuration yields a membrane dominated load-carrying behavior.

Due to the rectangular topology of the roof domain a single patch representation is quite natural. The surveyed quantity is the vertical deflection $\displacePComp_z$ at the midpoint of the free edges. In \cite{belytschko_1985} the value of the reference solution is reported as $0.3024$, this values is calculated using a very fine mesh.

In \figref{fig:roof_conv_Mmixed} the displacement convergence of the $\bendingMoment$-mixed formulation for $p=1,2,3,4$ is shown. As expected, it tuns out that the convergence becomes considerably faster with higher polynomial degree. For $p=3$ and $p=4$ already the fourth refinement step with $11$ control points per side yields a relative error of less than $1\%$, whereas the discretizations with $p=2$ and especially $p=1$ require a much finer mesh ($19$ and $400$ control points per side, respectively) to provide the same accuracy.

According to \tabref{tab:roof}, for $\vp=2$ the results obtained in \cite{echter_2013} for the standard purely displacement-based 3-parameter Kirchhoff-Love shell formulation (3p) conform well with the $\bendingMoment$-mixed shell elements developed in this paper for fine discretizations but provide worse results for coarse meshes.

\begin{table}
\centering
\begin{tabular}{p{4cm}p{1cm}p{1cm}p{1cm}p{1cm}p{1cm}p{1cm}}
\hline
Control points per edge & $5$ & $9$ & $13$ & $20$ & $25$ & $30$ \\
\hline
$p=1$\\
$\bendingMoment$-mixed &  $0.0252$ & $0.0471$ & $0.0665$ & $0.1027$ & $0.1291$ & $0.1528$ \\ 
\hline
$p=2$ \\ 
$\bendingMoment$-mixed &  $0.1028$ & $0.2636$ & $0.2940$ & $0.2997$ & $0.3003$ & $0.3004$ \\
3p (Echter \cite{echter_2013})    &  $0.0440$ & $0.2077$ & $0.2801$ & $0.2975$ & $0.2994$ & $0.3004$ \\
\hline
\end{tabular}
\caption{Scordelis-Lo roof, displacements ($\bendingMoment$-mixed, 3p).}
\label{tab:roof}
\end{table}

\begin{figure}[!ht]
	\centering
	\begin{subfigure}[b]{0.49\textwidth}
		\includegraphics[width=\textwidth]{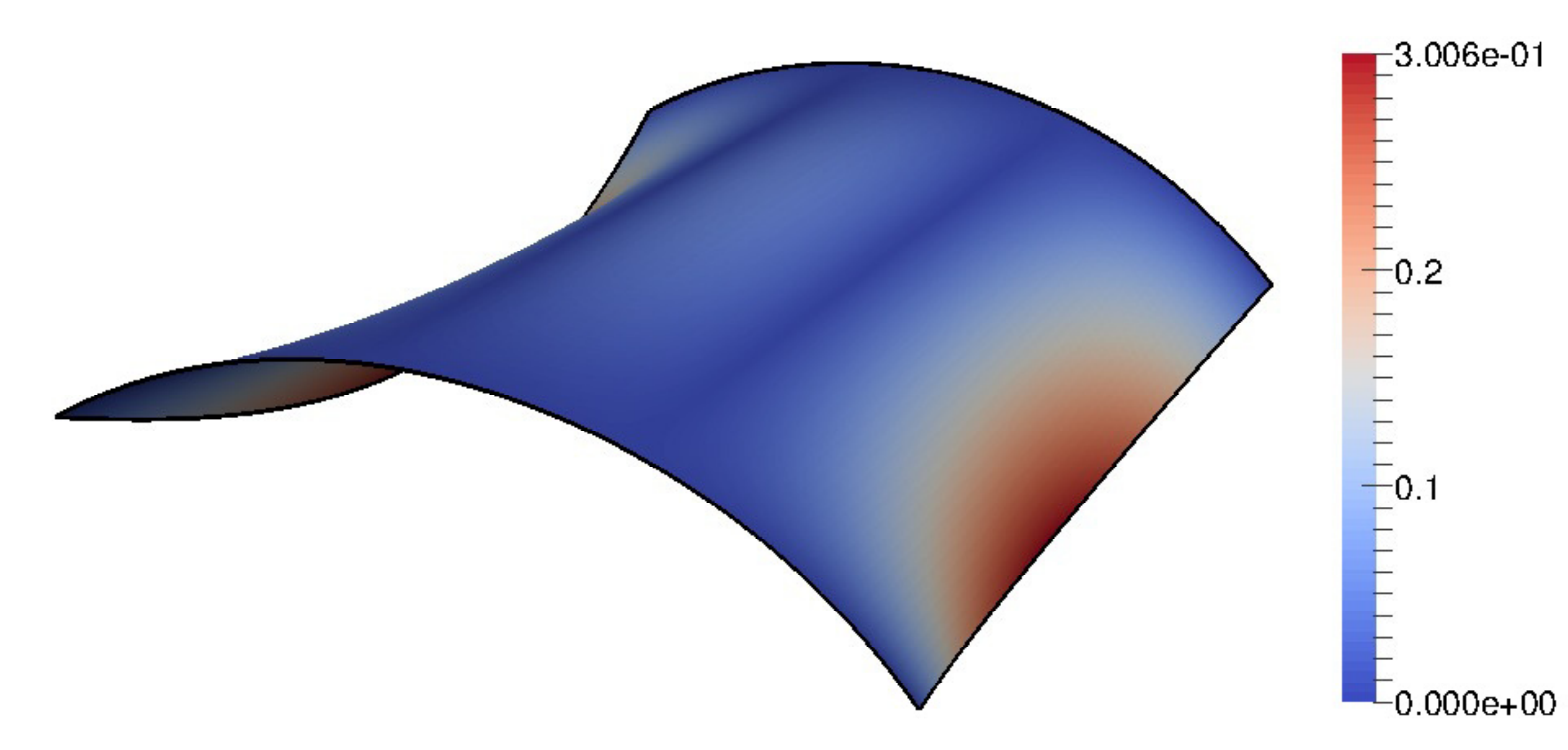}
		\includegraphics[width=\textwidth]{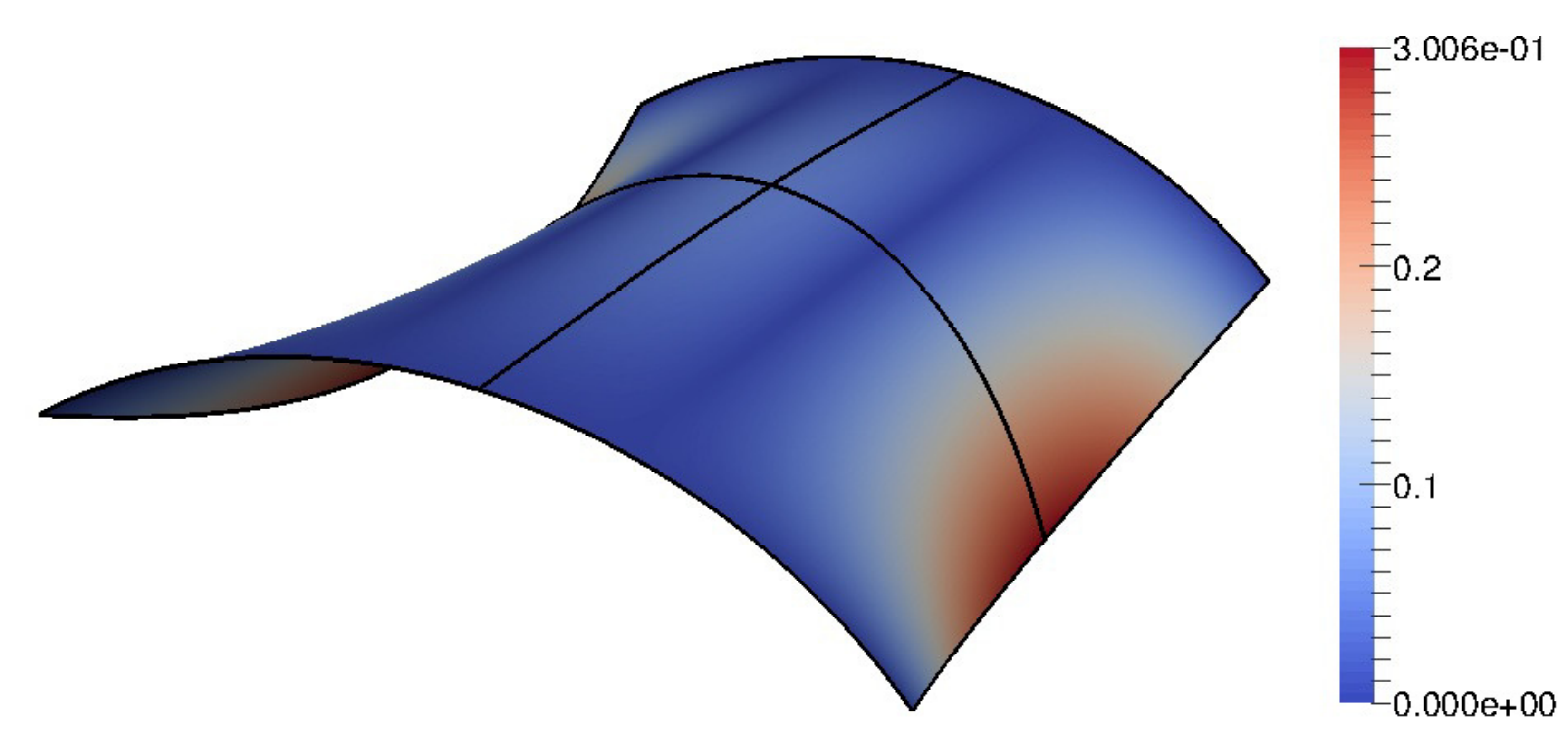}
		\caption{Displacement $\displacePComp_z$.}
	\end{subfigure}
	\begin{subfigure}[b]{0.49\textwidth}
		\includegraphics[width=\textwidth]{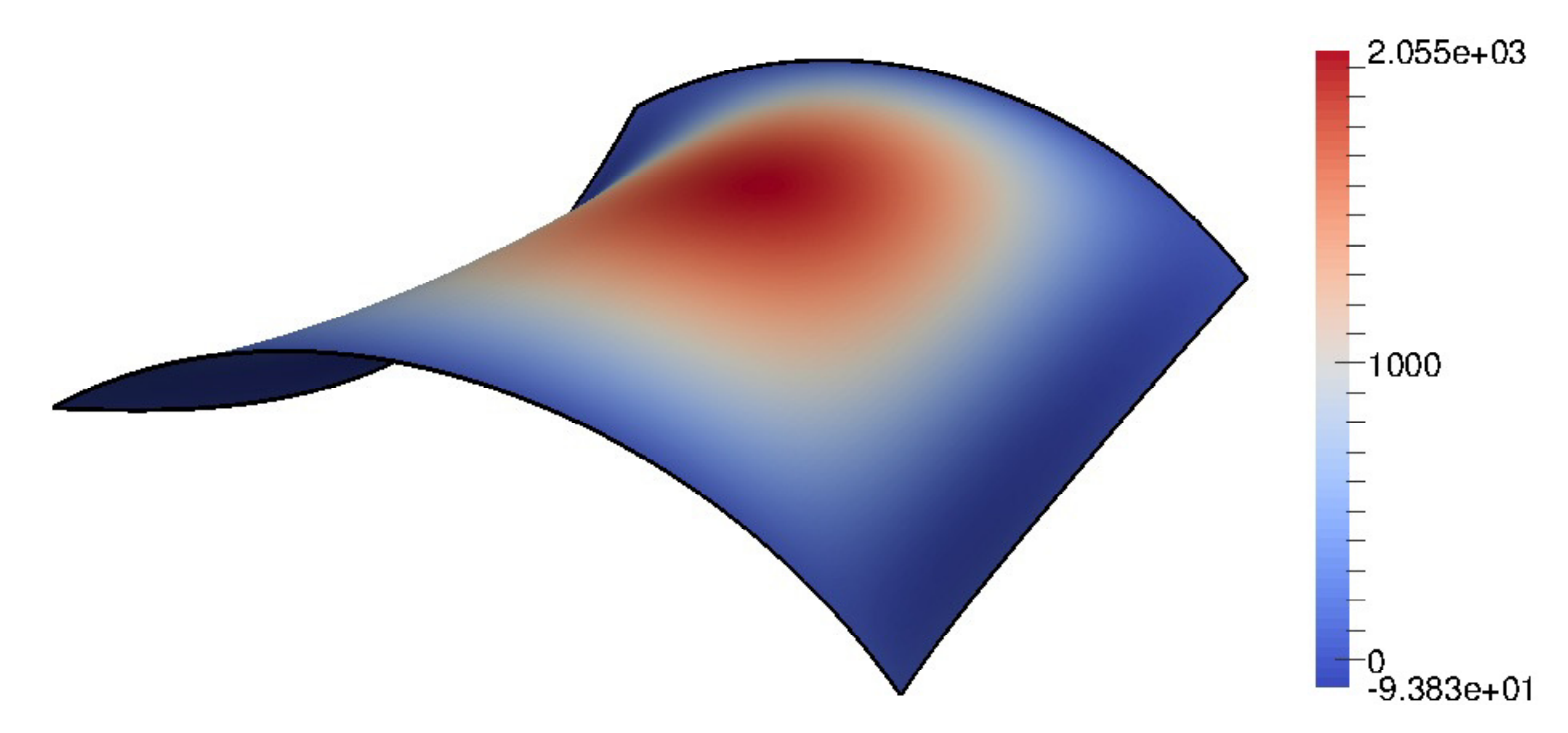}
		\includegraphics[width=\textwidth]{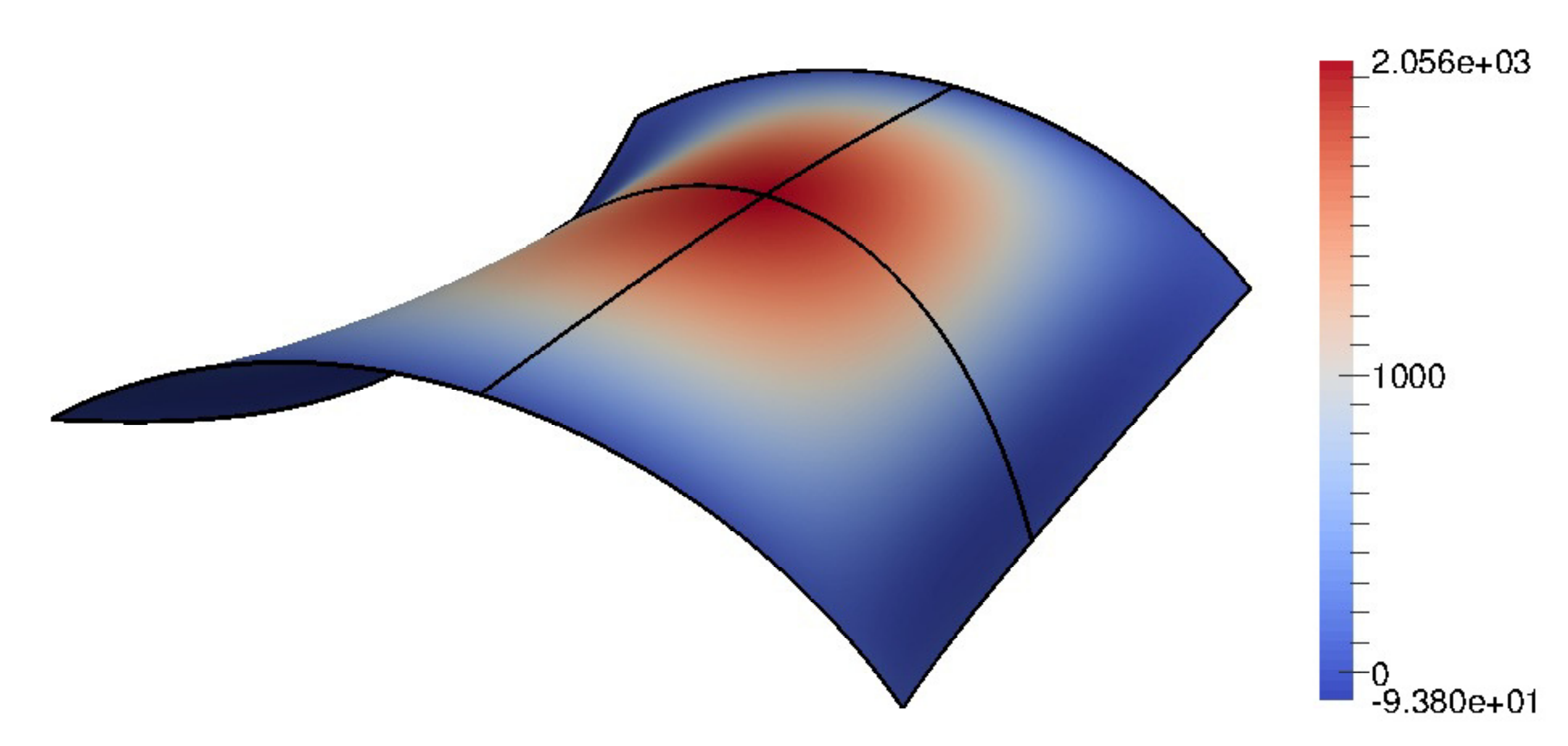}
		\caption{Bending moment component $\bendingMomentPComp^{xx}$.}
	\end{subfigure}
	\caption{Scordelis-Lo roof, analysis results for one patch (upper) and four patch (lower) geometry representations.}
	\label{fig:roof_patches}
\end{figure}

For comparison we consider a second representation of the midsurface using four patches, as illustrated in \figref{fig:roof_geometry}. In \figref{fig:roof_patches} vertical displacements and bending moments for the one patch and four patch geometry representation are compared. Visually no difference like discontinuity across patch interfaces can be seen. As in \cite{echter_oesterle_bischoff_2013}, in order to obtain $\bendingMomentPComp^{xx}$ the components of the bending moment tensor defined in the curvilinear coordinate system are first transformed into a local Cartesian basis with $\overline x$ and $\overline z$-axis aligned with $\xi^1$- and $\xi^3$-directions, for details see \cite{bischoff_wall_bletzinger_ramm_2004}.  
The obtained convergence results are quite similar to the ones received for the one patch geometry representation in \figref{fig:roof_conv_Mmixed}, and are therefore not shown.

\subsubsection{Pinched hemisphere}
\label{sec:hemisphere}
\begin{figure}[!ht]
	\centering
	\begin{subfigure}[b]{0.44\textwidth}
		\includegraphics[width=\textwidth]{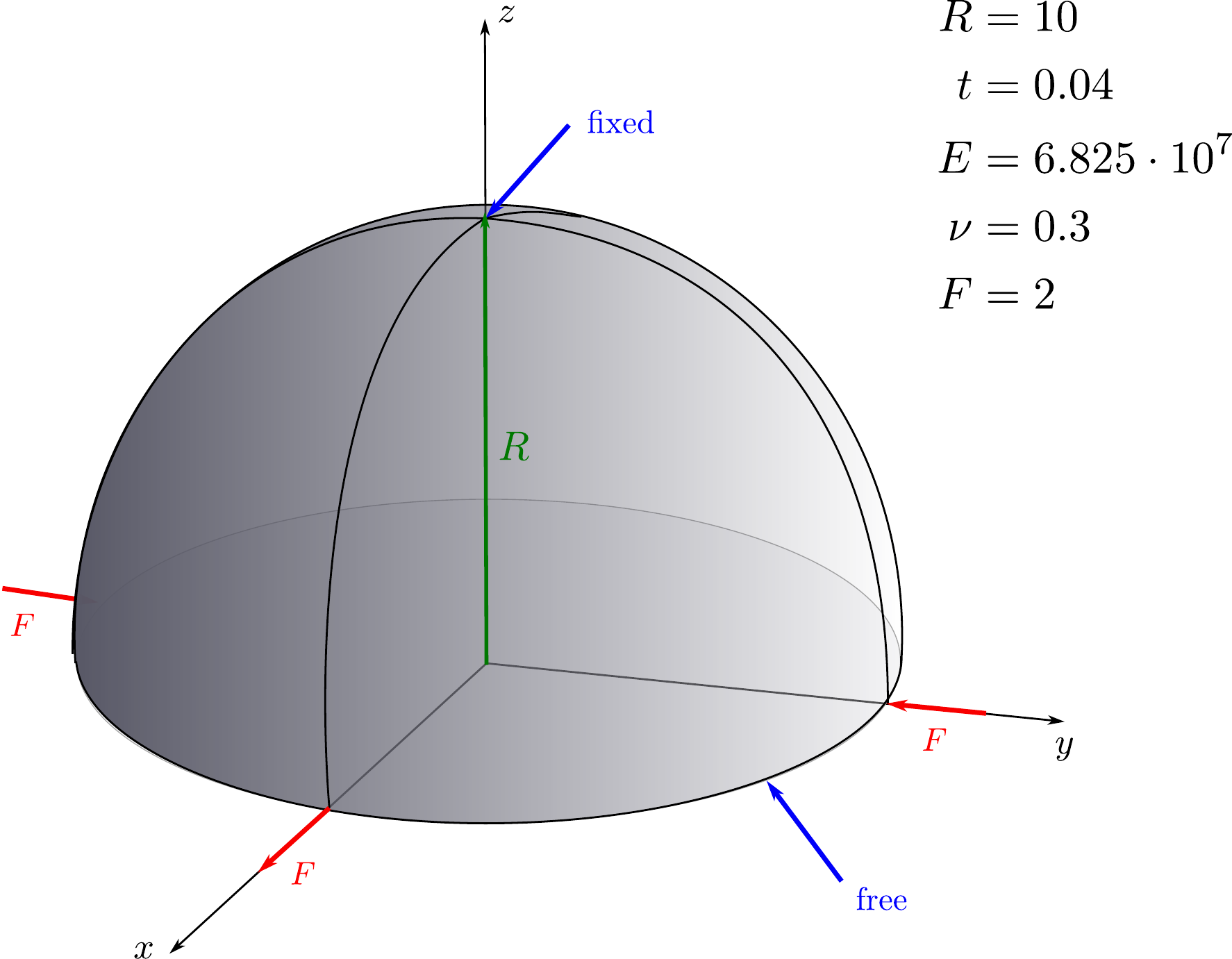}
		\caption{Problem setup.}
		\label{fig:hemisphere_geometry}
	\end{subfigure}
	\begin{subfigure}[b]{0.55\textwidth}
		\includegraphics[width=\textwidth]{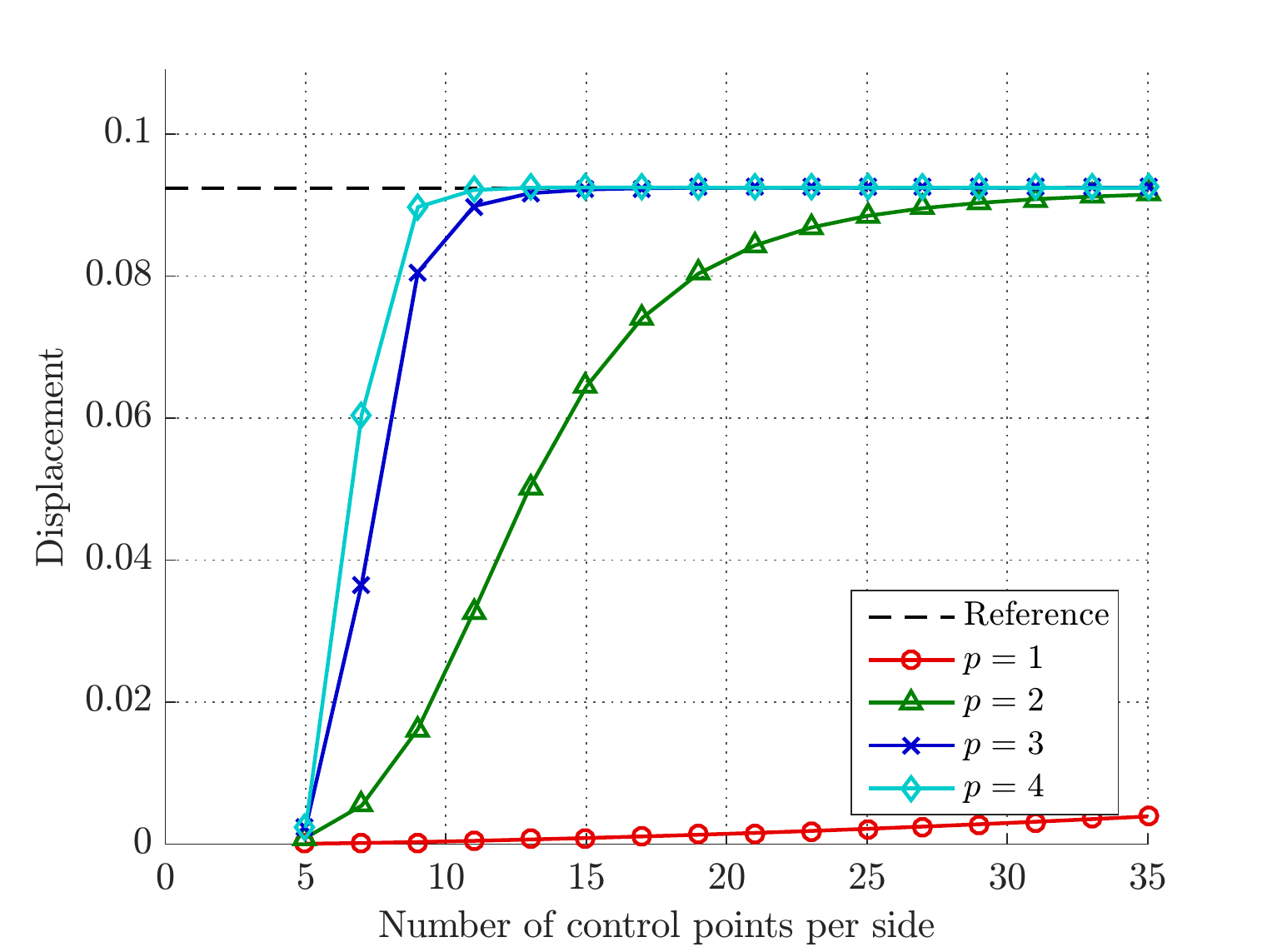}
		\caption{Displacement convergence of $\bendingMoment$-mixed formulation.}
		\label{fig:hemisphere_conv_Mmixed}
	\end{subfigure}
	\caption{Pinched hemisphere.}
\end{figure}

The problem setup of the pinched hemisphere benchmark is shown in \figref{fig:hemisphere_geometry}. The structure is fixed at the top and free along the bottom circumferential edge. The shell has a slenderness of $\frac{R}{t} = 250$, with radius $R$ and thickness $t$ as defined in \figref{fig:hemisphere_geometry}, and is subject to four radial point loads $F= \pm 2$ at its bottom. This configuration yields a bending dominated behavior with almost no membrane strains.

The undeformed midsurface is modeled using four patches, as illustrated in \figref{fig:hemisphere_geometry}. The investigated quantity is the radial displacement at the points where the loads are applied, with the value of the reference solution reported as $0.0924$ in \cite{belytschko_1985}.

In \figref{fig:hemisphere_conv_Mmixed} the displacement convergence of the $\bendingMoment$-mixed formulation for $p=1,2,3,4$ is shown. Most observations carry over from the Scordelis-Lo roof benchmark to the pinched hemisphere. The slowed convergence of low order discretizations with $\vp=1$ and $\vp=2$ becomes even more amplified. The discretization with $p=2$ requires $35$ control points per edge (for each of the four patches) to reach a relative error of less than $1\%$. The discretization with $p=1$ shows a considerably slower convergence. With a reasonable number of control points no acceptable accuracy could be achieved. For $p=3$ and $p=4$ only $13$ control points per edge are needed to obtain the same accuracy. 

The reason for the slow convergence in case of $p=1$ and $p=2$ for the Scordelis-Lo roof and (even more severe) for the pinched hemisphere is membrane locking, which is mechanically the inability to represent pure bending without unwanted, parasitic membrane strains, see, e.g., \cite{bischoff_wall_bletzinger_ramm_2004} for more details on membrane locking. Similar results showing evidence of membrane locking have already been observed in \cite{kiendl_2009, echter_2013}.
We will come back to the important issue of membrane locking in \sectref{sec:membraneLocking}.

\subsubsection{Pinched cylinder}
\label{sec:cylinder}
\begin{figure}[!ht]
	\centering
	\begin{subfigure}[b]{0.44\textwidth}
		\includegraphics[width=\textwidth]{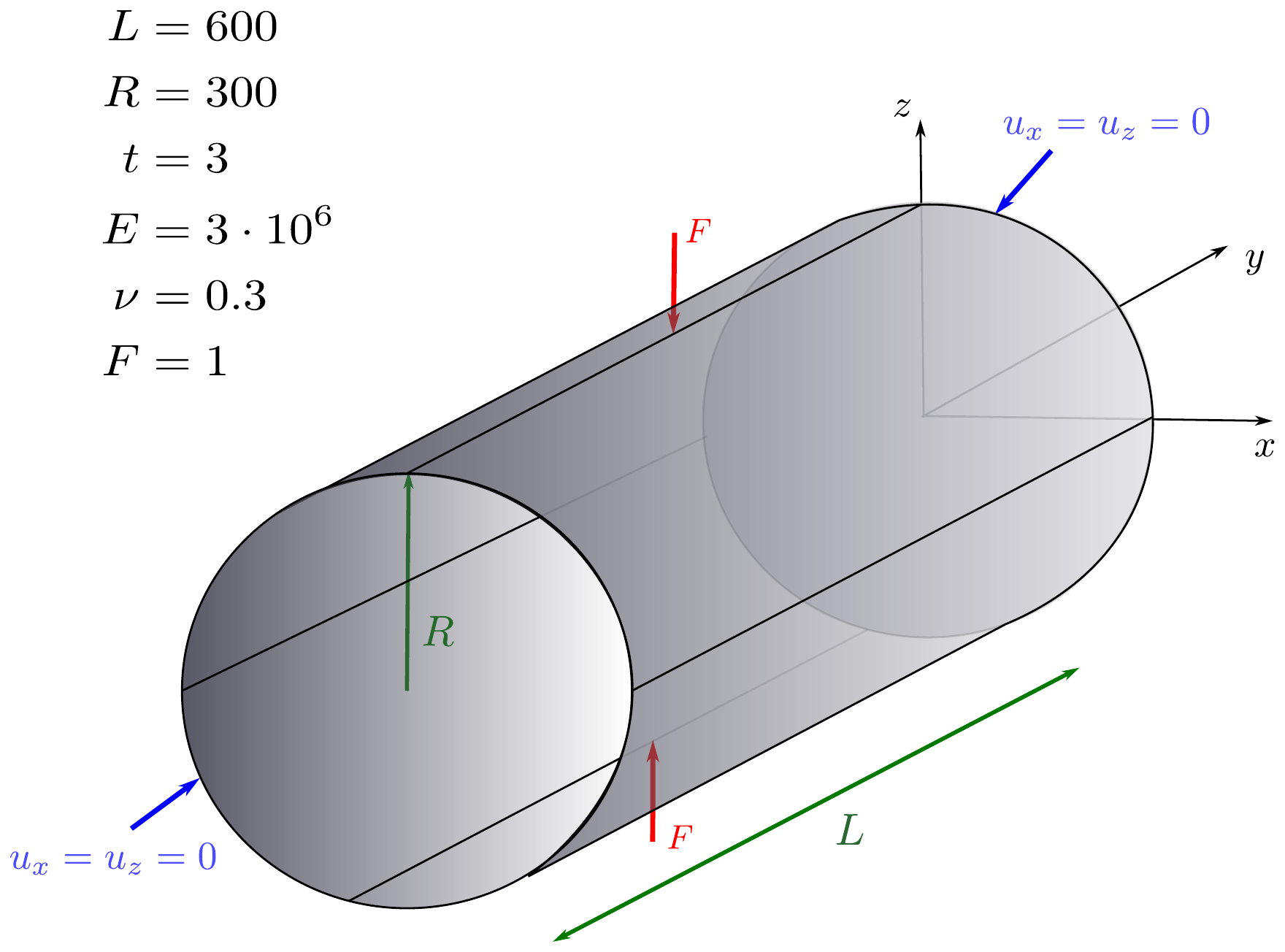}
		\caption{Problem setup.}
		\label{fig:cylinder_geometry}
	\end{subfigure}
	\begin{subfigure}[b]{0.55\textwidth}
		\includegraphics[width=\textwidth]{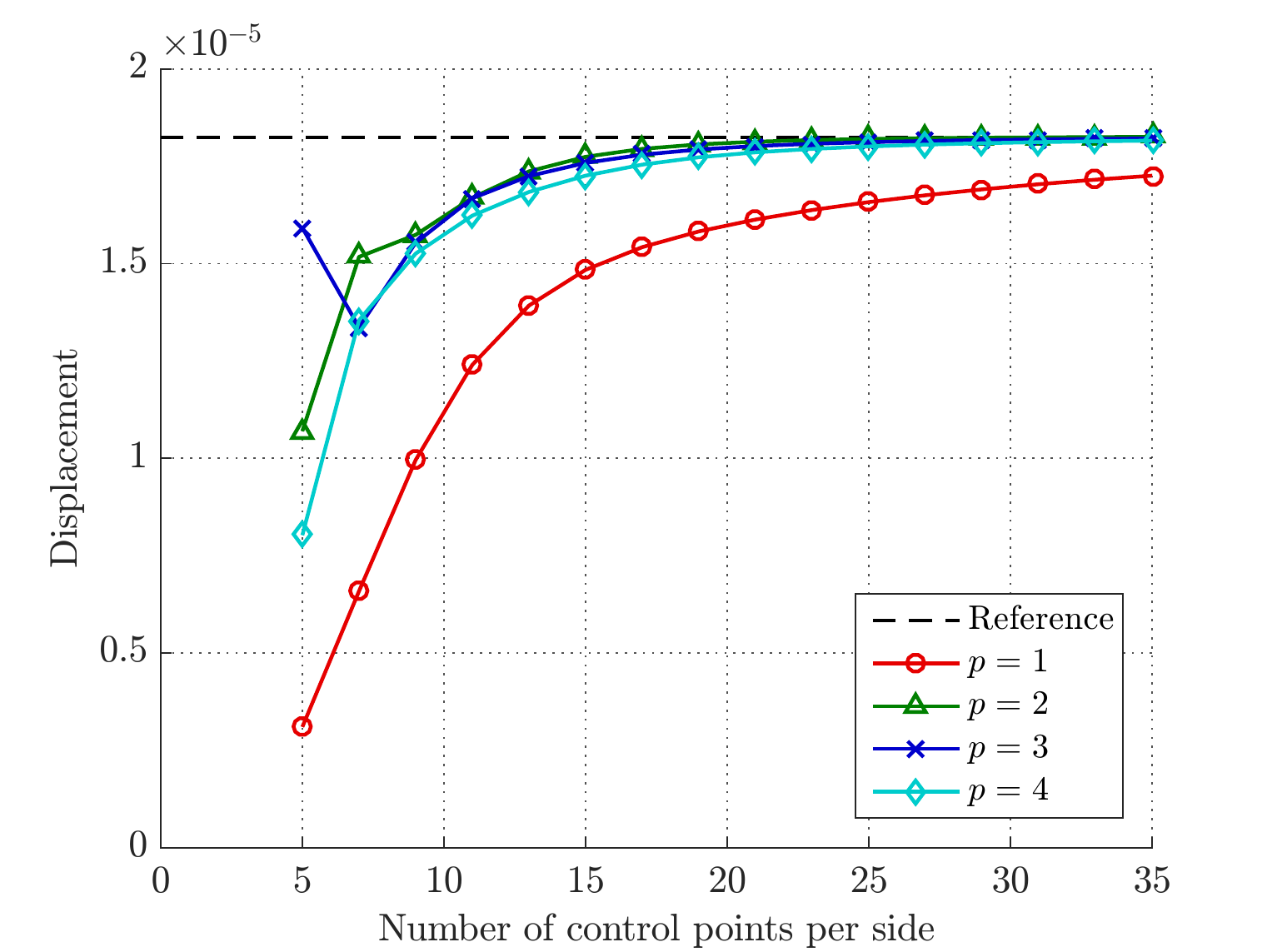}
		\caption{Displacement convergence of $\bendingMoment$-mixed formulation.}
		\label{fig:cylinder_conv_Mmixed}
	\end{subfigure}
	\caption{Pinched cylinder.}
\end{figure}

The last benchmark of the shell obstacle course is the pinched cylinder. The problem setup is shown in \figref{fig:cylinder_geometry}. The structure is supported with rigid diaphragms at both ends. The cylinder has moderate slenderness $\frac{R}{t} = 100$, with radius $R$ and thickness $t$ as defined in \figref{fig:cylinder_geometry}, and is subject to two opposite point loads $F= \pm 1$ in the middle. This configuration yields a severe test for both inextensional bending and complex membrane states. 

The undeformed midsurface is modeled using four patches, as illustrated in \figref{fig:cylinder_geometry}. The analyzed quantity is the radial displacement at the points where the loads are applied, with the value of the reference solution reported as $1.8248 \cdot 10^{-5}$ in \cite{belytschko_1985}.

In \figref{fig:cylinder_conv_Mmixed} the displacement convergence of the $\bendingMoment$-mixed formulation for $p=1,2,3,4$ is shown. The convergence behavior differs from the ones obtained for the first two benchmarks. The differences in the results due to the usage of different polynomial orders are significantly smaller, especially discretizations with $p=2,3,4$ provide comparable results. Note, already the first-order discretization shows an acceptable performance.   

\subsection{Membrane locking and \texorpdfstring{$\bendingMoment$-$\membraneForce$}{MN}-mixed formulation}
\label{sec:membraneLocking}
\begin{figure}[!ht]
	\centering
	\begin{subfigure}[b]{0.44\textwidth}
		\includegraphics[width=\textwidth]{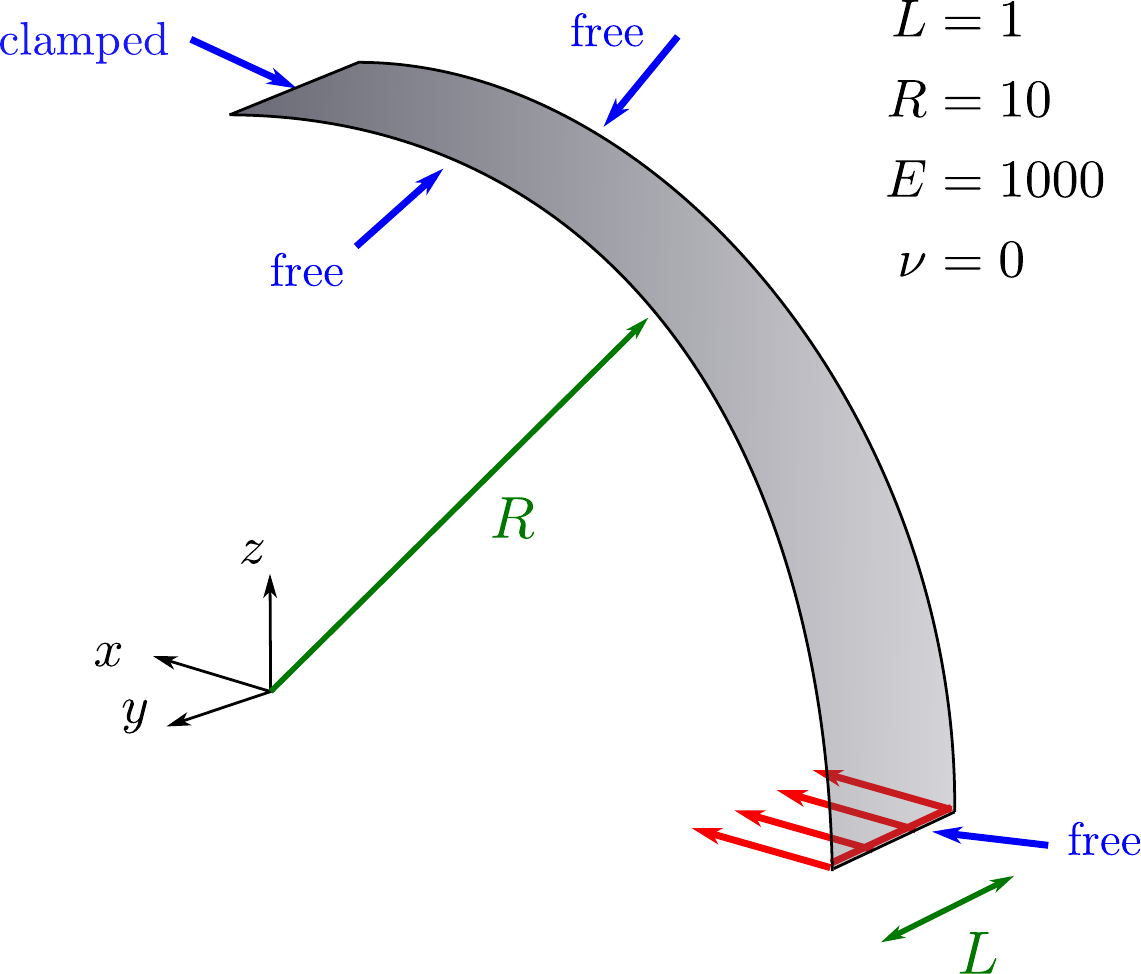}
		\caption{Problem setup.}
		\label{fig:cylinderStrip_geometry}
	\end{subfigure}
	\begin{subfigure}[b]{0.55\textwidth}
		\includegraphics[width=\textwidth]{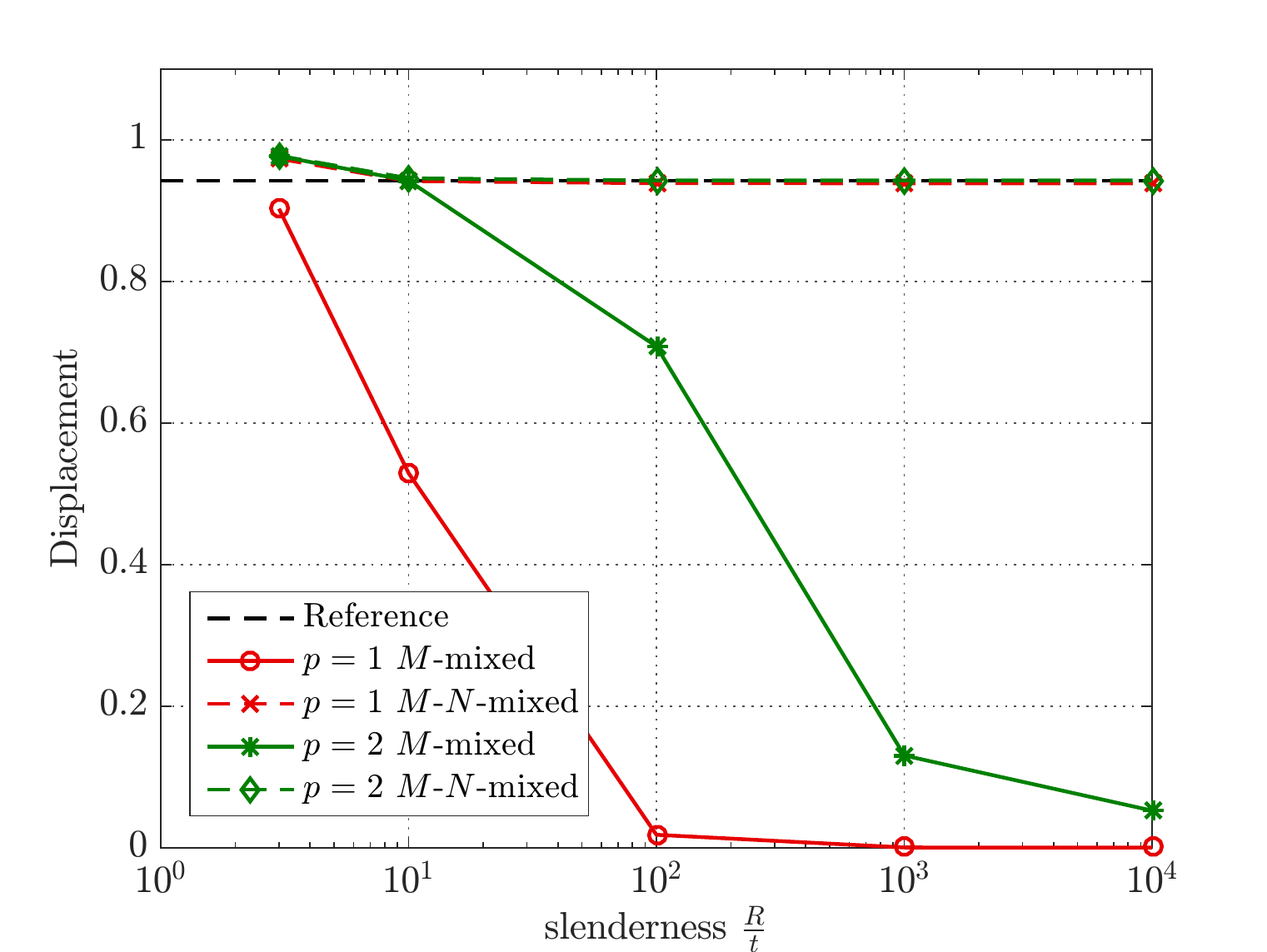}
		\caption{Displacement convergence.}
		\label{fig:cylinderStrip_locking}
	\end{subfigure}
	\caption{Cylindrical shell strip.}
\end{figure}

In order to investigate membrane locking we consider first a simple model problem consisting of a cylindrical shell strip, see, e.g.,\cite{bischoff_wall_bletzinger_ramm_2004, echter_oesterle_bischoff_2013}. The problem setup of this example is shown in \figref{fig:cylinderStrip_geometry}. The structure is clamped along the edge $x=0$ and subject to a constant constant line load in radial direction with magnitude $q_x = 0.1 \cdot t^3$ at the opposite free edge.  This configuration yields a bending dominated behavior. Therefore, membrane locking has to be expected in case the applied discrete formulation is not free from membrane locking.

The quantity of interest is the radial displacement at the midpoint of the free edge, with exact value $0.942$ independent of the slenderness $\frac{R}{t}$, according to an analytical solution based on Bernoulli beam theory (cf. \cite{echter_oesterle_bischoff_2013}). In order to ensure comparability with \cite{echter_oesterle_bischoff_2013} the domain is discretized with a mesh of 10 elements in longitudinal direction and one element in the other direction.

In \figref{fig:cylinderStrip_locking} the influence of varying slenderness $\frac{R}{t}$ on the obtained displacement is investigated for the original $\bendingMoment$-mixed formulation and the extended $\bendingMoment$-$\membraneForce$-mixed formulation for $p=1$ and $p=2$. Both results obtained with the $\bendingMoment$-mixed formulation show severe membrane locking. The radial displacement tends to zero as the slenderness $\frac{R}{t}$ is increased. In case of $p=2$ already for a moderate slenderness of $\frac{R}{t} = 100$ unphysical membrane strains lead to a considerable underestimation of the tip displacement of approximately $25\%$ and for $p=1$ the behavior is even worse. However, for both polynomial orders the results obtained with the $\bendingMoment$-$\membraneForce$-mixed formulation indicates that this formulation completely removes the undesired membrane locking effects.

According to \tabref{tab:cylinderStrip}, for $p=2$ the results in \cite{echter_oesterle_bischoff_2013} for the standard purely displacement-based 3-parameter formulation (3p) and the corresponding formulation with a mixed Hybrid Stress modification of the membrane part (3p-HS) conform well with the $\bendingMoment$-mixed and $\bendingMoment$-$\membraneForce$-mixed shell elements proposed in this paper. 

\begin{table}
\centering
\begin{tabular}{p{4.5cm}p{1cm}p{1cm}p{1cm}p{1cm}p{1.5cm}}
\hline
Slenderness $\frac{R}{t}$ & $10$ & $100$ & $1000$ & $10 000$ \\
\hline
$p=2$\\ 
$\bendingMoment$-mixed  & $0.9420$ & $0.7075$ & $0.1156$ & $0.0112$ \\ 
$\bendingMoment$-$\membraneForce$-mixed & $0.9454$ & $0.9423$ & $0.9422$ & $0.9422$ \\ 
$3$p (Echter et al. \cite{echter_oesterle_bischoff_2013})   & $0.9326$ & $0.6635$ & $0.0225$ & $0.0002$ \\ 
$3$p-HS (Echter et al. \cite{echter_oesterle_bischoff_2013})  & $0.9386$ & $0.9425$ & $0.9425$ & $0.9425$ \\ 
\hline
\end{tabular}
\caption{Cylindrical shell strip, displacements ($\bendingMoment$-mixed, $\bendingMoment$-$\membraneForce$-mixed, 3p, 3p-HS). }
\label{tab:cylinderStrip}
\end{table}

\begin{figure}[!ht]
	\centering
	\begin{subfigure}[b]{0.49\textwidth}
		\includegraphics[width=\textwidth]{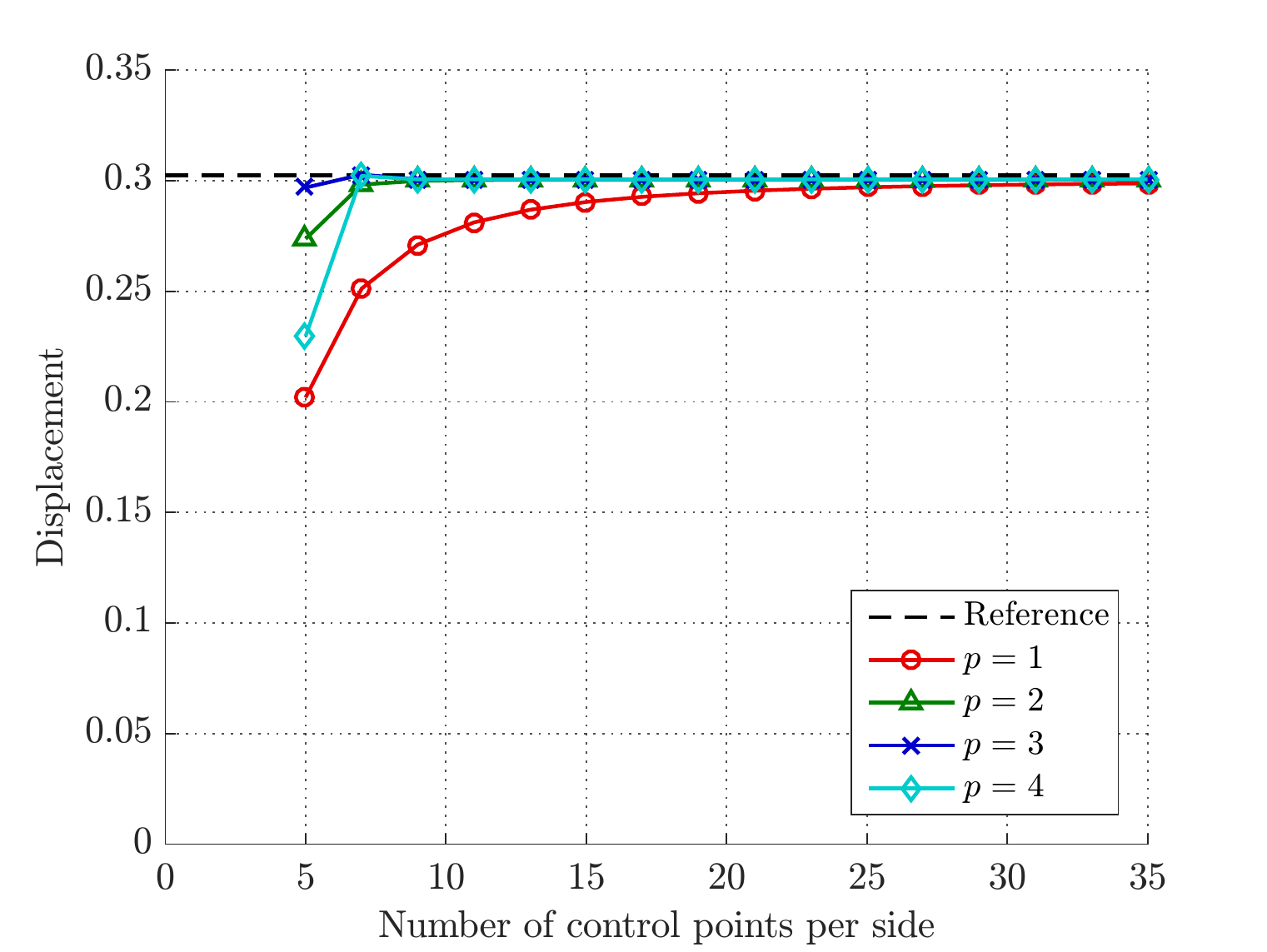}
		\caption{Scordelis-Lo roof}
		\label{fig:cylinder_conv_MNmixed}
	\end{subfigure}
	\begin{subfigure}[b]{0.49\textwidth}
		\includegraphics[width=\textwidth]{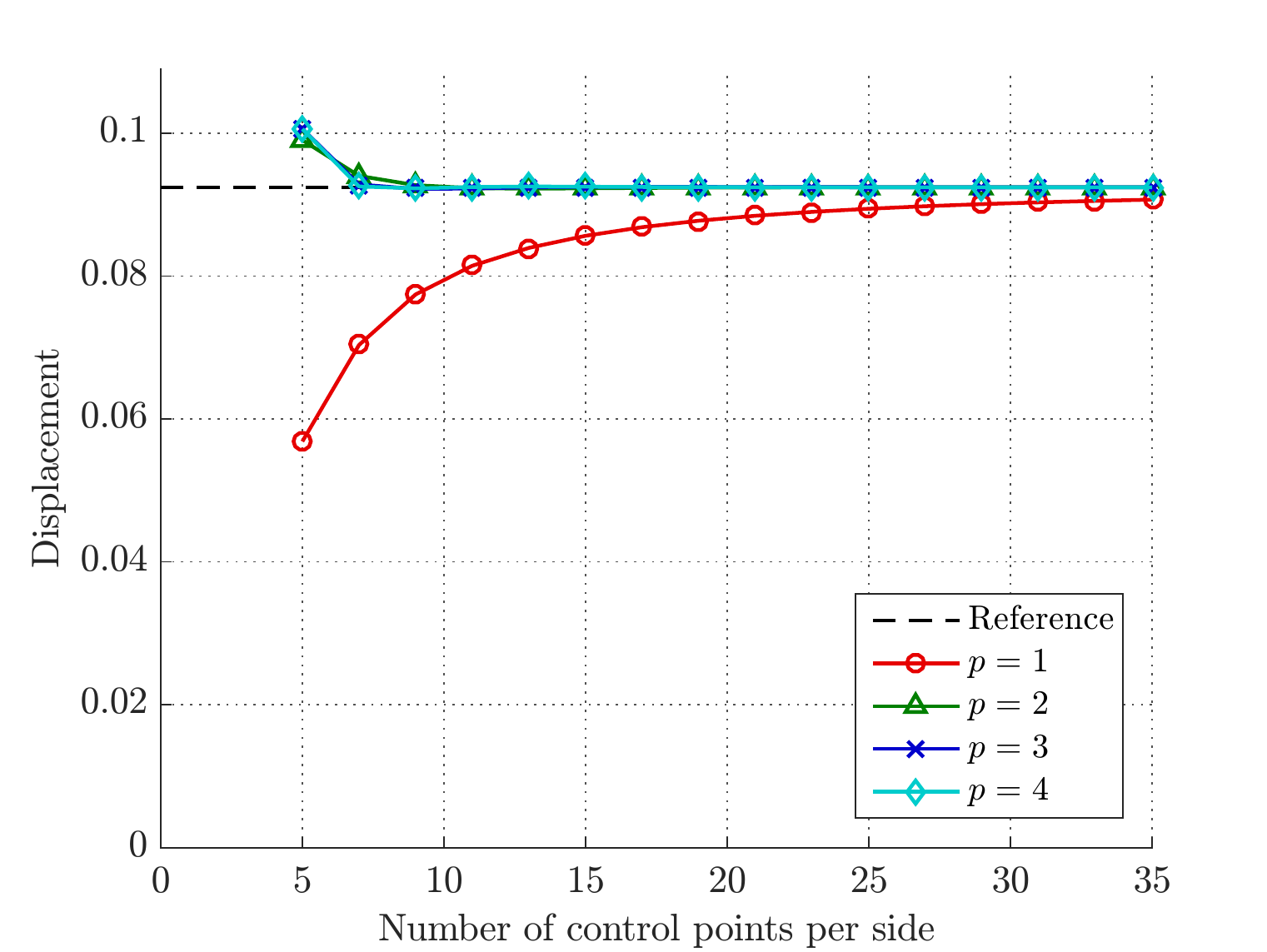}
		\caption{Pinched hemisphere}
		\label{fig:hemisphere_conv_MNmixed}
	\end{subfigure}
	\caption{Displacement convergence of $\bendingMoment$-$\membraneForce$-mixed formulation.}
\end{figure}

The improved convergence behavior of the $\bendingMoment$-$\membraneForce$-mixed formulation observed in the basic cylindrical shell strip problem carries over to the Scordelis-Lo roof and pinched hemisphere benchmark (cf. \sectref{sec:roof} and \sectref{sec:hemisphere}). In \figref{fig:cylinder_conv_MNmixed} and \figref{fig:hemisphere_conv_MNmixed} the displacement convergence of the $\bendingMoment$-$\membraneForce$-mixed formulation for $p=1,2,3,4$ is shown for the Scordelis-Lo roof and pinched hemisphere benchmark, respectively. For both problems the locking free $\bendingMoment$-$\membraneForce$-mixed shell element converges considerably faster to the reference value for the same number of control points compared to the $\bendingMoment$-mixed element in \figref{fig:roof_conv_Mmixed} and \figref{fig:hemisphere_conv_Mmixed}. Furthermore, the differences in the results for different polynomial orders are significantly reduced.

In \tabref{tab:Nmixed} the results of the $\bendingMoment$-$\membraneForce$-mixed formulation are compared to the results for the 3p-HS formulation of Echter \cite{echter_2013}. Recall, our formulation does not require an $H^2$-conforming discretization in contrast to the one in \cite{echter_2013}. A quantification of the improvements in numbers is obtained by comparing \tabref{tab:Nmixed} with \tabref{tab:roof}.


\begin{table}
\centering
\begin{tabular}{p{4cm}p{1cm}p{1cm}p{1cm}p{1cm}p{1cm}p{1cm}}
\hline
Control points per edge & $5$ & $9$ & $13$ & $20$ & $25$ & $30$ \\
\hline
$p=1$\\
$\bendingMoment$-$\membraneForce$-mixed &  $0.2020$ & $0.2711$ & $0.2869$ & $0.2941$ & $0.2970$ & $0.2978$ \\
\hline
$p=2$ \\ 
$\bendingMoment$-$\membraneForce$-mixed &  $0.2737$ & $0.2999$ & $0.3005$ & $0.3006$ & $0.3006$ & $0.3006$ \\
3p-HS (Echter \cite{echter_2013})  &  $0.2517$ & $0.3000$ & $0.3005$ & $0.3006$ & $0.3006$ & $0.3006$ \\
\hline
\end{tabular}
\caption{Scordelis-Lo roof, displacements ($\bendingMoment$-$\membraneForce$-mixed, 3p-HS).}
\label{tab:Nmixed}
\end{table}

\section{Concluding remarks and future work}

The numerical experiments in \sectref{sec:numerical_experiments} demonstrate that the proposed $\bendingMoment$-mixed shell element works in well-known benchmark problems. The results for the $\bendingMoment$-$\membraneForce$-mixed formulation indicate that this formulation is free from membrane locking.

In contrast to Kirchhoff-Love type thin shells, for conforming Reissner-Mindlin shell elements standard $C^0$-continuous shape functions are commonly used. In order to avoid transverse shear locking, in \cite{long_bornemann_cirak_2012,echter_oesterle_bischoff_2013} an $H^2$ based Mindlin-Reissner formulation obtained by a change of variables from rotations to transverse shear strains is proposed. This approach, which was later mathematically analyzed in the context of plates in \cite{lovadina_2015}, has the advantage that transverse shear locking is eliminated on the continuous formulation level, independent of a particular discretization. An application of our technique to obtain an $H^1$-mixed-formulation to this $H^2$ based Mindlin-Reissner formulation seems possible and would be worth to investigate.

\bibliographystyle{abbrv}
\bibliography{KLshell_bibliography}

\end{document}